\documentclass[oneside,draft,11pt]{amsart}
\usepackage{amssymb, amscd}
\usepackage[all]{xy}
\usepackage{verbatim}

\usepackage{latexsym}
\usepackage{amssymb}

\title[]{Central values of L-functions over CM fields}

\address{Department of Mathematics\\
University of Michigan\\
2074 East Hall\\
Ann Arbor, MI 48109}
\email{hxue@umich.edu}

\usepackage[dvips]{epsfig}
\usepackage[all]{xy}
\usepackage{amsmath,amsopn,latexsym,amssymb,amscd,amsthm}

\author{Hui Xue}
\newtheorem{theorem}{Theorem}[subsection]
\newtheorem{prop}[theorem]{Proposition}
\newtheorem{lem}[theorem]{Lemma}
\newtheorem{cor}[theorem]{Corollary}

\newtheorem{remark}{Remark}[subsection]

\numberwithin{equation}{subsection}
\newtheorem{thm}{Theorem}

\newcommand{\tr}{{\mbox{tr}}}
\newcommand{\N}{{\bf{N}}}
\newcommand{\BZ}{{\mathbb{Z}}}

\newcommand{\BQ}{{\mathbb{Q}}}

\newcommand{\BC}{{\mathbb{C}}}
\newcommand{\BA}{{\mathbb{A}}}
\newcommand{\BR}{{\mathbb{R}}}
\newcommand{\BH}{{\mathbb{H}}}

\newcommand{\D}{\displaystyle}
\begin{document}

\maketitle
\tableofcontents

\begin{abstract}
In the paper we prove an explicit formula  for the central values of certain Rankin L-functions. These L-functions are L-functions of Hilbert newforms over a totally real field F, twisted by unitary  Hecke characters of a totally imaginary quadratic extension of F. The formula generalizes our former result on L-functions twisted by finite characters. It also eliminates the restriction on ramifications of Popa's paper.     
\end{abstract}

\section{Introduction} \label{s1}
In this paper we will prove an explicit formula for central values of certain Rankin $L$-functions.  The setting is as follows. Let $f$ be a Hilbert newform over a totally real field $F$ with trivial central character and has weight $(2k_1,2k_2,\cdots,2k_d)$. Suppose the level of $f$ is $N$. Let $\chi$  be a unitary character on $\BA^{\times}_K/K^{\times}\BA_F^{\times}$ of conductor $c(\chi)$. Here $K$ is a totally imaginary quadratic extension of $F$,  $\BA_F^{\times}$ and $\BA^{\times}_K$ are the idele groups of $F$ and $K$ respectively. Let $\pi$ and $\pi_{\chi}$ be the automorphic representations of $GL_2(\BA_F)$ associated to $f$ and $\chi$ respectively, then we use $L(s,\pi\times\pi_{\chi})$ to denote the Rankin-Selberg convolution of $L(s,\pi)$ and $L(s,\pi_{\chi})$. The $L$-function $L(s,\pi\times\pi_{\chi})$ has a functional equation of the form $L(s,\pi\times\pi_{\chi})=\epsilon(s,\pi\times\pi_{\chi}) L(1-s,\pi\times\pi_{\chi})$. 

{\bf Assumption.}  We  assume that $2$, $N$ and $d_{K/F}$ (the absolute relative discriminant of $K/F$) are co-prime to each other.

Under this assumption the sign $\epsilon(1/2,\pi_f\times\pi_{\chi})$ is given by $(-1)^{|\Sigma|}$, and $\Sigma$ is the following set of places of $F$
\begin{equation*}
\Sigma=\Sigma_1\cup
\{\text{$v|N$ such that $\omega_v(N)=-1$}\},
\end{equation*}
where $\Sigma_1$ is a set of archimedean places determined by  weights of $\chi$ and $f$ (see the following paragraph), $\omega$ is the quadratic character of $F^{\times}\backslash \BA^{\times}$ associated to $K/F$. In this paper we will find an explicit formula for $L(1/2,\pi\times\pi_{\chi})$ when the sign of the functional equation is $+1$, namely when $|\Sigma|$ is even.

To determine $\Sigma_1$ we assume that the  character $\chi$ has the infinity type 
\begin{align} \label{s1,1}
\chi_{\infty}=\sum_{\tau\in\Phi}r_{\tau}\tau+\sum_{\rho\tau}r_{\rho\tau}\rho\tau,
\end{align}
where $\Phi$ is a fixed CM type of $K$, $\rho$ is the complex conjugation. Since $\chi$ is unitary we have and assume $-r_{\tau}=r_{\rho\tau}>0$. Let $r_{v}=r_{\rho\tau}$ if $v$ is the infinite place of $F$ under $\tau$, then $\Sigma_1$ is the set of $v$'s such that $r_v<k_v$. The rest subset of infinite places of $F$ is denoted by $\Sigma_2$.

If $|\Sigma|$ is even  there is a unique (up-to isomorphism) quaternion algebra $B$ over  $F$ such that $B$ is ramified exactly over $\Sigma$. The CM field $K$ can be embedded in $B$ and the embedding will be fixed from now on. It is well-known that 
there is exactly one irreducible automorphic representation $\pi^B$ of $B^{\times}(\BA)/\BA^{\times}$, such that  $\pi^B$ and $\pi$ are related  by  Jacquet-Langlands correspondence.

We now state a short version of  our main theorem. For the complete statement and notations see Theorem \ref{thm} in Section \ref{s4}.

\begin{thm} 
 The central value of $L(s,\pi\times\pi_{\chi})$ is given by 
$$L(1/2,\pi\times\pi_{\chi})= C\cdot \frac{(\varphi^{\ast},\varphi^{\ast})}{(\varphi^B,\varphi^B)}\cdot\left|\int_{K^{\times}\BA_F^{\times}\backslash\BA^{\times}_K}\varphi^B(t)\chi(t)^{-1}dt\right|^2,$$
where $C$ is an explicitly determined positive constant, $\varphi^{\ast}$ is a certain automorphic form in the representation space of $\pi$, and $\varphi^B$ is a certain automorphic form in the representation space of $\pi^B$. 
\end{thm}

Our  paper is devoted to the proof of this theorem.  When $\chi$ is a finite character we have obtained such a formula  in \cite{Xue02} by using a geometric argument following \cite{Zhang01} (initiated by Gross \cite{Gross87}).  In \cite{Popa03}, whose approach will be taken in our paper, a similar formula over real quadartic fields is proved under the assumption that $N$ is squarefree and $\chi$ is unramified. The goal of our paper is to obtain a similar result over CM fields and remove the restriction on conductors. 

We would  like to thank Shouwu Zhang for his  encouragement and help.

\subsection{Notations} \label{s1,s1}
We  use $F$ to denote either a totally real number field or a non-archimedean local field, depending on the context. When $F$ is a non-archimedean local field, we write $\varpi$, $\mathcal{O}$, and $U$ to denote its fixed uniformizer, ring of integers and group of units respectively. The letter $K$  denotes either a totally imaginary quadratic extension of $F$ when $F$ is a number field, or a quadratic extension (split or not) of $F$ when $F$ is a local field. If $S$ is  an algebraic object defined over $F$, for instance if $S$ is a vector field or an algebraic group  over $F$, then we use $S_{\BA}$ to denote its adelic points over the adele ring $\BA=\BA_F$.  The space of Schwartz functions on a topological space $V$ is denoted by $S(V)$.  For a subset $A$ of $V$ its characteristic function is denoted by $1_A$.

The notation $U_0(1)$ represents the standard maximal compact subgroup $GL_2(\hat{\mathcal{O}}_F)$ of $GL_2(\BA)$, while $U_0(N)$ denotes the matrices in $U_0(1)$ with lower left entries divisible by $N$. Similar notations hold in local situation. We denote by $Z$, $B$ and  $N$ the center, the standard Borel and unipotent subgroup of $GL_2$ respectively. The notion $n(x)$ denotes a unipotent matrix with $x$ in the upper right. Sometimes we also denote $GL_2$ by $G$. We let $T_1$ be the diagonal subgroup of $GL_2$ with lower right entry equal to $1$, and let $t(a)$ be the matrix in $T_1$ with $a$ in the upper left. 

We fix an additive character $\psi$ of $\BA/F$ by the formula $\psi(a)=\psi_{\BQ}(\tr_{F/\BQ}(a))$, here $\psi_{\BQ}$ is the standard additive character of $\BA_{\BQ}$ and $\tr_{F/\BQ}$ denote the trace map from $\BA/F$ to $\BA_{\BQ}/\BQ$. Therefore the conductor $\delta$  of $\psi$ is the relative discriminant of $F/\BQ$. Locally we also fix an unramified additive character  $\psi^0$ of $F_v$. The character $\psi^0$ and $\psi_v$ (the restriction of $\psi$ on $F_v$) are related by the identity $\psi(x)=\psi^0(\delta_v x)$, where $\delta_v$ is the discriminant of $F_v$ over $\BQ_v$.

We will use $\N$ and $\tr$ (with or without the subscript $V$) to denote the reduced norm and trace of elements in a quadratic or quaternion algebra $V/F$.
 We use $\omega_V$ to denote the quadratic character associated to a quadratic space $V$ and use $\omega$  especially to denote the quadratic character of the quadratic extension $K/F$. Finally we use $c(\pi)$ (or $c(\omega)$) to denote the conductor of the representation $\pi$ (or the character $\omega$).

\section{Preliminaries} \label{s2}

In this section we will first recall some facts on Weil representations of $SL_2$ and theta kernels. Then we will construct some special local Schwartz functions associated to these representations. 
\subsection{Weil representations and theta kernels} \label{s2,s1}
First let $F$ be a local field of characteristic $0$. Let $\psi$ be a fixed nontrivial additive character of $F$. Let $V$ be a $2n$-dimensional space over $F$ with a non-degenerate quadratic form $q$.  
Let $GO(V)$ be the group of similitudes of the quadratic space $V$ defined by
$$GO(V)=\{\sigma\in GL_F(V):\,q(\sigma v)=\nu(\sigma)q(v)\,\text{with $\nu(\sigma)\in F^{\times}$ for $v\in V$}\},$$
where $\nu$ is called the similitude factor.
We write 
$$R(GO(V), GL_2)=\{(h,g)\in GO(V)\times GL_2: \nu(h)=\det(g)\}.$$
Note that $R(GO(V), GL_2)$ is isomorphic to the semidirect product $GO(V)\ltimes SL_2$ by the map
$$(h,g)\mapsto \left(h,\begin{pmatrix}1&0\\0&\det(g)^{-1}\end{pmatrix}\cdot g\right).$$
Since ${\rm dim}(V)$ is even, one can construct a Weil representation $r_{q}$ defines a representation of $SL_2(F)$  on the space $S(V)$ of Bruhat-Schwartz functions on $V$ by the following rule
$$r_{q}\left(\begin{pmatrix}1& a\\ 0&1\end{pmatrix}\right)f(x)=\psi(aq(x))f(x),$$
$$r_{q}\left(\begin{pmatrix}a&0\\ 0& a^{-1}\end{pmatrix}\right)f(x)=|a|^n\omega_V(a)f(ax),$$
$$r_{q}\left(\begin{pmatrix}0&1\\ -1& 0\end{pmatrix}\right)f(x)=\gamma\hat{f}(x),$$
here $\gamma$ is an eighth root of unity, $\omega_V$ is the quadratic character of $F^{\times}$ associated to the quadratic space $V$. The measure on $V$  is taken to be self-dual with respect to the Fourier transform $\hat{f}(x)=\int_Vf(y)\psi(q(x,y))dy$, where $q(x,y)=q(x+y)-q(x)-q(y)$.

If $r'=r_{\lambda}$ is the Weil represnetation attached to $(V,\lambda q)$ then
$$\omega'=\omega,\,\,\gamma'=\omega(\lambda)\gamma,\,\, dx'=|\lambda|dx.$$

If $r_q'$ is the Weil representation attached to $(V,q)$ but with respect to $\psi'(x)=\psi(\delta x)$, then
\begin{equation} \label{s2,s1,3}
r'(g)=r(t(\delta)gt(\delta)^{-1}).
\end{equation}

The representation $r_{q}$ can be extended to a representation $r_q$ of $R(GO(V),GL_2)$ on the same space $S(V)$ by letting 
$$r_q(h,g)f(x)=L(h)r_q(g_1)f(x),$$
where $\D g_1=\begin{pmatrix} 1&0\\ 0& \det(g)^{-1}\end{pmatrix}\cdot g\in SL_2(F)$ and $L(h)f(x)=|\nu(h)|^{-n}f(h^{-1} x)$.

In this paper we only consider two types of $V$. 
The first type is $V=K$ for a quadratic algebra $K$ over $F$ (split or not), equipped with a quadratic form $q(x)=\Lambda \N_{K}$ for certain $\Lambda\in F$.  
The group $K^{\times}$ is a subgroup of index $2$ of $GO(V)$ whose action on $K=V$ is given  by multiplication. The  similitude factor   is given by  $\nu(t)=\N_{K}(t)$.  The whole group $GO(V)$ is the semi-direct product of $K^{\times}$ and the involution group of $K/F$. 
The other type is $(V,q)=(B,\N_B)$, where $B$ is a quaternion  algebra over $F$. The similitude group $GO(V)$ is the semi-direct product of its  connected identity component $GSO(V)$ and the standard involution of $B$. Also, one has $GSO(V)=B^{\times}\times B^{\times}/F^{\times}$ through the map $(g_1,g_2)\mapsto g_1b g_2^{-1}$ for $b\in B$. Under this identification the similitude factor on $GSO(B)$ is given by $\N_B(g_1g_2^{-1})$.

Now  suppose $F$ is totally real number field. Let $(V,q)$ be an anisotropic quadratic space of dimension $2n$ over $F$.  A similar global construction as the above gives the Weil representation of the group $R(GO(V_{\BA}),GL_2(\BA))$ for a fixed addtive character  $\psi$ of $\BA/F$.

The theta kernel on $R(GO(V_{\BA}),GL_2(\BA))$  is defined by
\begin{equation} \label{s2,s1,2}
\theta(h,g;\phi)=\sum_{x\in V_F}r_q(h,g)\phi(x)
\end{equation}
for any $\phi\in S(V_{\BA})$. 

There are three global quadratic spaces to be considered. The first one is  $(K,\N_K)$ for the CM field $K$, and it corresponds to the construction of theta series associated to $\chi$, see Section \ref{s4,s1}. Anothe one is $(K,\Lambda \N_K)$  and corresponds to the construction of an Eisenstein series (Section \ref{s4,s1}). The third one is $(B,\N_B)$ for the quaternion algebra $B$ over $F$ and corresponds to the theta lift between $GSO(B)$ and $GL_2$ (Section \ref{s4,s3}).

\subsection{Special vectors} \label{s2,s3} 
For later applications we need to find some special vectors $\phi$ in the space $S(V_{\BA})$ for various quadratic spaces $V$. We choose and fix a global embedding $K\to B$ such that $B=K+Kj$, where $j$ is an  element in $B$ with $j^2=-\Lambda\in \mathcal{O}_F$ and $j^{\iota}=-j$ ($\iota$ is an involution of $B$ with respect to the fixed embedding of $K$). As quadratic spaces over $F$ one has an orthogonal decomposition
$(B,\N_B)=(K,\N_K)\oplus(K,\Lambda\N_K)$. 

 Let $K=F+Fi$, where  $i^{\iota}=-i$ and $i^2=n_K$. So the quaternion algebra  $B$ has a decomposition $B=F\oplus Fi\oplus Fj\oplus Fk$, where $i^2=n_K$, $j^2=-\Lambda$, and $ij=-ji=k$. By the assumption we have $\epsilon(B_v)=(n_K,-\Lambda)_v$, where $(\cdot,\cdot)_v$ is the local Hilbert symbol of $K$ at $v$.

\begin{lem} \label{s2,s3,l1}
Let $c={\rm Max}(v(c(\pi)),v(c(\pi_{\chi}))={\rm Max}(v(N),v(D))$ with $D=c(\chi)^2c(\omega)$, then
we can choose $j\in B$ such that $\Lambda=\N(j)$ falls in one of  the following five cases for $m\in\BZ_{+}$:\\
(1)  $c=0$, with $v$  inert in $K$, then $v(\Lambda)=2m$;\\
(2)  $c=0$, with $v$ split in $K$, then $v(\Lambda)=m$;\\
(3a)  $c=2t+1$, with $v$ inert in $K$, then $v(\Lambda)=1$;\\
(3b)  $c=2t\ge2$, with $v$ inert in $K$, then $v(\Lambda)=0$;\\
(4)  $c\ge1$, $v$ is split in $K$, then $v(\Lambda)=m=0,1$;\\
(5)  $c=1$, $v$ is ramified in $K$, then $v(\Lambda)=0$. 
\end{lem}
\begin{proof}
We first show that one can remove all the even positive prime powers in $\Lambda$ at places where $c$ is nonzero. Let $\frak{p}$ be a prime such that $\frak{p}^{2t}|\Lambda$ for a maximal $t\ge1$. By Chebotarev density theorem one can find an odd prime $\frak{q}$ (which is coprime to every ramified place) such that $\frak{p}{\frak{q}}^{-1}=x\mathcal{O}_F$ with $x\in F$. If we let $j'=j/x^{t}$ and $\Lambda'=-(j')^2$, then  the new $\Lambda'$ has the power less than $2$ at the prime  $\frak{p}$. Obviously $K+Kj'$ equals $B$. 

If $K/F$ is inert at a place $v$ (or equivalently, at a prime $\frak{p}$), then  $\epsilon(B_v)=(n_K,\Lambda)_v=(-1)^{v(\Lambda)}$. But by our construction we have  $(-1)^c=\epsilon(B_v)$, so $v(\Lambda)$ and $c$ have the same parity.So we have got  the first 4 cases. 

Now we assume that $c=1$ and $K_v/F_v$ is ramified. By the above argument we may also assume $v(\Lambda)=0,1$. By Chebotarev we find a good prime  $\frak{Q}$ of $K$ and an $x\in K$ such that $x\mathcal{O}_K=\frak{Q}\frak{P}^{-1}$, where $\frak{P}$ is the prime ideal of $K$ which divides $\frak{p}$ (the prime of $v$).  If we take $j'=x^{-1}j$ and $\Lambda'=-j'^2=-\N_K(x^{-1})j^2$, then $B=K+Kj'$ and $v(\Lambda')=0$.      
\end{proof}
    
Note that the case $v(N)\ne 0$ corresponds to Cases 3 and 4, and the case $v(c(\chi))\ne0$ corresponds to Cases 3b and 4.

The theta series and Eisenstein series to be studied later on depend on the choice  of $\phi\in S(V_{\BA})$. 
There are three cases. The
 first one
 is the theta lift from  $ GSO(K)\backslash GSO(\BA_K)$ to $G(F)\backslash G(\BA)$. The second one occurs in the definition of the Eisenstein series $f(s,g;\phi)$. The third one is the theta lift  from $G(F)\backslash G(\BA)$ to $ GSO(B)\backslash GSO(B_{\BA})$ (see Section \ref{s4}).

{\bf Nonarchimedean case.}
Now we construct special Schwartz functions at a finite place $v$ according to the above three cases. Let $K=K_v$ and $F=F_v$. 

{\bf N1.} Let $V=K$ be the fixed quadratic extension of $F=F_v$ equipped with $\N_{K/F}$ as the quadratic form. If $v\not| c(\chi)$ then  we take $\phi_1(x)=1_{\mathcal{O}_K}(x)$, if $v|c(\chi)$ then we take $\phi_1(x)=\chi(x)1_{\mathcal{O}_K^{\times}}$. 

{\bf N2.} 
 We let  $\phi_2\in S(K)$  be the characteristic function $1_L$ of a lattice $L$ constructed as follows (according to the five cases discussed in Lemma \ref{s2,s3,l1}):\\
(1):  $L={\varpi^{-m}\mathcal{O}_K},$\\
(2):  $L={\mathcal{O}}\oplus{\varpi^{-m}\mathcal{O}}$,\\
(3a) and (3b): $L={\varpi^t\mathcal{O}_K},$\\ 
(4):  $L={\mathcal{O}}\oplus{\varpi^{c-m}\mathcal{O}}$, where $v(\Lambda)=m$ for $m=0,1$ \\
(5):  $L={\mathcal{O}_K},$ where $v(\Lambda)=0$. 

Actually in Case 2, if we take $\varpi_K=(1,\varpi)$ be a uniformizer of $K=F\oplus F$ then $L=\varpi^{-m}\mathcal{O}_K$, so Case 2 and Case 1 can be treated together. The  remark also applies to Case 4, which can be treated similarly as Case 3.

Let $R=\mathcal{O}_{K}+Lj$, a simple calculation shows  the following. 

\begin{lem} \label{s2,s3,l2}
The lattice $R$ is an (Eichler) order of reduced discriminant $\varpi^c$ in $B$.
\end{lem}

{\bf N3.} Let  $V=B$ be a quaternion algebra with its norm $\N_B$ as the quadratic form. Under the orthogonal decomposition $B=K+Kj$ we take $\phi(x_1+x_2j)=\phi_1(x_1)\phi_2(x_2)$, where $\phi_1$ and $\phi_2$ are constructed in {\bf N1} and {\bf N2} respectively. So the function $\phi$ is given by
\begin{equation*}
\phi(x)=\begin{cases}1_{R} & \text{if $v\not|c(\chi)$,}\\
 \chi(x_1)1_{R^{\times}} & \text{if $v|c(\chi)$ and $x=x_1+x_2j$.}
\end{cases}
\end{equation*}
The last line is valid at $v|c(\chi)$ as we can write $Lj=\mathcal{O}_Kj'$ such that $j'^2$ is a unit multiple of $c(\chi)^2$. In the seond case we also denote $\phi$ by $\chi$, which is a character on $R^{\times}$.

 If $K/F$ is ramified  we let $\tilde{R}$ denote the following maximal order of $M_2(F_v)$
\begin{align} \label{s2,s3,1}
\tilde{R}=\{ a+bj:a,b\in\delta_{K/F}^{-1}\text{ such that }a-bu\in\mathcal{O}_{K_v}\},
\end{align} 
where $\delta$ is the different of $K/F$, and $u\in U_{K}$ is such that $j^2=\N_K u$. Here such an $u$ exists because $(n_K,j^2)_v=(n_K,-\Lambda)_v=1$, i.e, $j^2$ is a norm in $K$. We write $\phi'=1_{\tilde{R}}$.

Now back to the global setting. Let $\psi$ be a fixed additive character of $\BA/F$ of conductor $\delta$. We form a compact subgroup of $B^{\times}_{\BA_f}$
$$\widehat{R}^{\times}=\prod_{v\not|c(\omega)}R_v^{\times}\prod_{v|c(\omega)}\tilde{R}^{\times}_v,$$
and define a character $\chi$ (by abuse of language) on $\widehat{R}^{\times}$ by
$\chi(k)=\prod_{v|c(\chi)}\chi_v(k_v),$
where $R_v$, $\widetilde{R}_v$ and $\chi_v$ are defined as above.
 For $x=\prod_{v\not|\infty}x_v\in B_{\BA_f}$, we define 
$$\phi'_f(x)=\prod_{v\not|c(\omega)}\phi_v(x_v)\cdot\prod_{v|c(\omega)}\phi_v'(x_v).$$
 Now we determine the level structure of $\phi'_f$  in the Weil representation  $r_B$ attached to $(B,\N_B)$. 
\begin{prop}[Level structure] \label{s2,s3,p1}  For $\kappa\in t(\delta)^{-1}U_0(N c(\chi)^2)t(\delta)$ and $k_1,k_2\in \widehat{R}^{\times}$ such that $\det \kappa=\N_B(k_1 k_2^{-1})$
$$r_B((k_1,k_2),\kappa)\phi'_f=\chi(k_1^{-1}k_2)\phi'_f,$$
where $(k_1,k_2)$ is regarded as an element in $GSO(B_{\BA_f})=\BA_f^{\times}\backslash B^{\times}_{\BA_f}\times B^{\times}_{\BA_f}$.
\end{prop}
\begin{proof}
As $r_B((k_1,k_2),\kappa)\phi'(x)=r_B(\kappa_1)\phi'(k_1^{-1}x k_2)$ and $\phi'(k_1^{-1}xk_2)=\chi(k_1^{-1}k_2)\phi'$, we are reduced to check
\begin{equation} \label{s2,s3,5}
r_B(\kappa)\phi'=\phi',
\end{equation}
 for any $\kappa\in t(\delta)^{-1}U_0(N c(\chi)^2)t(\delta)\cap SL_2(\BA_f)$.
 It suffices to check (\ref{s2,s3,5}) locally.  At a finite place $v$,  by (\ref{s2,s1,3}) we may assume $\psi_v$ is unramified and check (\ref{s2,s3,5}) for generators of  $U_0(N c(\chi)^2)$. 

If $v(Nc(\chi)^2)=0$, then $\phi_v$ (or $\phi'_v$) is the characteristic function for a maximal $\mathcal{O}_v$-order of $B_v$ and the check of (\ref{s2,s3,5}) is simple (or see \cite{Watson02}). 

We assume now $v(Nc(\chi)^2)=c>0$. 
\begin{lem} 
Let $r_{\lambda}$ as usual be the Weil representation associated to $(K_v,\lambda\N)$. Assume  $K_v/F_v$ is unramified (split or not).  If $\phi_v=1_{\varpi_{K_v}^{c-M}\mathcal{O}_{K_v}}$ and $v(\lambda)=M$, then $r_{\lambda}(\kappa)\phi_v=\phi_v$ for any $\kappa\in U_0(\varpi^{c})\cap SL_2(F_v)$.

If $\phi_v=\chi(x)1_{\mathcal{O}_K^{\times}}$ for $v(c(\chi))\ne0$, then $r_1(\kappa)\phi_v=\phi_v$ for any $\kappa\in U_0(c(\chi)^2)\cap SL_2(F_v)$.
\end{lem}
\begin{proof}
The subscript $v$ will be suppressed during the proof. The group $U_0(\varpi^{c})\cap SL_2(F)$ is generated by matrices of the type
$$l(a)=\begin{pmatrix}a&0\\0&a^{-1}\end{pmatrix},\,\, n(b)=\begin{pmatrix}1&b\\0&1\end{pmatrix},\,\, m(u\varpi^c)=\begin{pmatrix}1&0\\u\varpi^c&0\end{pmatrix},$$
where $a, u\in \mathcal{O}^{\times}$, $b\in\mathcal{O}$. It is clear that $r_{\lambda}(l(a))\phi=\phi$ and $r_{\lambda}(n(b))\phi=\phi$. For $m(u\varpi^c)$ we note that
$$m(u\varpi^c)=-w n(-u\varpi^c) w,\,\, {\rm for}\,\, w=\begin{pmatrix}0&1\\-1&0\end{pmatrix}.$$
Now $r_{\lambda}(w)\phi=\gamma\hat{\phi}$, here $\hat{\phi}$ is the Fourier transform with respect to the self-dual measure $dx$. So $\hat{\phi}={\rm const}\cdot 1_{\varpi^{-M}\overline{\varpi}_K^{-c+M}}$, where  $\overline{\varpi}_K$ is the conjugate of $\varpi_K$ over $F$. It is easy to see $r_{\lambda}(n(-u\varpi^c))$ fixes $\hat{\phi}$, hence ($\gamma=1$ in our case)
$$r_{\lambda}(m(u\varpi^c))\phi(x)= r_{\lambda}(w)\hat{\phi}(x)=\phi(-x)=\phi(x)$$
by the Fourier inversion formula.

See \cite{Popa03} Proposition 2.5.1 for the proof of the second statement. 
\end{proof}
Consider the following decomposition of $r_B$ via  $(B_v,\N_{B_v})=(K_v,\N_{K_v})\oplus(K_v,\Lambda\N_{K_v})$
$$r_{B_v}(\kappa)\phi(x_1+x_2j)=r_1(\kappa)\phi_{1,v}(x_1)r_{\Lambda}(\kappa)\phi_{2,v}(x_2),$$
for $\kappa\in U_0(\varpi^c)\cap SL_2(F_v)$.
 The proof of Proposition \ref{s2,s3,p1} is  completed by  applying the  lemma to both $r_1(\kappa)$ and $r_{\Lambda}(\kappa)$.
 \end{proof}

{\bf Archimedean case.} \label{s2,s2}
Assume that we have fixed the global decomposition of $B=K+Kj$ induced by the chosen embedding of $K$ into $B$, such that $\N_B(j)=\Lambda$.
At an archimedean place $v$ there are two cases about the sign of $\Lambda$. If $k_v>r_v$, $B_v$ is definite, then $\Lambda_v$ has to be positive; if $k_v\le r_v$, $B_v$ is split and  $\Lambda_v$ is negative.
If we let $j_v=\begin{pmatrix}0&\sqrt{\Lambda_v}\\ -\sqrt{\Lambda_v}&0\end{pmatrix}$ then we have the identification $\BH\cong K_v\oplus K_vj_v$ by
$$\begin{pmatrix}u&v\\-\bar{v}&\bar{u}\end{pmatrix}\mapsto\begin{pmatrix}u&0\\ 0&\bar{u}\end{pmatrix}+\begin{pmatrix}v/\sqrt{\Lambda_v}&0\\ 0&\bar{v}/\sqrt{\Lambda_v}\end{pmatrix}j_v.$$
If $B_v\cong M_2(\BR)$, and $j_v=\begin{pmatrix}0&\sqrt{|\Lambda_v|}\\ \sqrt{|\Lambda_v|}&0\end{pmatrix}$, the identification becomes
$$\begin{pmatrix}u&v\\\bar{v}&\bar{u}\end{pmatrix}\mapsto\begin{pmatrix}u&0\\ 0&\bar{u}\end{pmatrix}+\begin{pmatrix}v/\sqrt{|\Lambda_v|}&0\\ 0&\bar{v}/\sqrt{|\Lambda_v|}\end{pmatrix}j_v.$$
In the following we will construct $\phi$'s on the spaces in the left side of the above identifications. In other words, we will construct  Schwartz functions on the quadratic spaces $\BC$ with the standard norm $\N$ (or $-\N$, see Case A2),  or $\BH$ and $M_2(\BR)$ with the standard reduced norms. See Remark \ref{s2,s2,r1} for the general case.

  Now let the additive character $\psi$ be defined by $\D \psi_v(x)=e^{2\pi ix}$.  

{\bf  A1.} Let $V=\BC$ be the quadratic space equipped with $q(z)=x^2+y^2$ over $F_v=\BR$, where $z=x_iy\in V$. Suppose the
character $ \chi(z)$ is given by $\D\frac{\bar{z}^r}{z^r}$, where $r$ is a
non-negative integer. Let us define 
$\phi_1(z)=2\bar{z}^{2r}e^{-2\pi|z|^2}\in{S}(\BC,\chi)$, where $S(\BC,\chi)$ represents the $\chi$-isotypic subspace of the whole Schwartz space.  Let $\D
X_+=\begin{pmatrix}0&1\\0&0\end{pmatrix}$ and $\D
X_-=\begin{pmatrix}0&0\\1&0\end{pmatrix}$ be two vectors in the Lie algebra $sl_2$ of $SL_2(\BR)$. In order to check the weights of $\phi$'s in the Weil representation we  need to check the action of $X_{+}-X_{-}$ on these vectors.

The effects of these two vectors on $S(\BC)$ are well-known (see \cite{Watson02}):
$$r(X_+)\phi=2\pi i(z\cdot \bar{z})\phi,$$ 
$$r(X_-)\phi=-\frac{1}{2\pi i}(\frac{\partial}{\partial
  z}\cdot\frac{\partial}{\partial {\bar{z}}})\phi.$$
It is easy to check that 
$$r(X_+-X_-)\phi_1=i(2r+1)\phi_1,$$
therefore $\phi_1(z)=2\bar{z}^{2r}e^{-2\pi|z|^2}$ or the corresponding theta lift has weight $2r+1$.

{\bf A2.} For the Eisenstein series, there are two cases. The first case is when the quadratic space is $\BC$ with the standard norm $q$.  Then we  take
$\D\phi_2(z)=p_l(4\pi |z|^2)e^{-2\pi|z|^2}\in S(\BC,1)$, where $l=k-r-1$ (when $k>r$) and $p_l(t)$ is the Laguerre polynomial of degree $l$:
$$p_l(t)=\sum_{j=0}^{l}\begin{pmatrix}l\\
  j\end{pmatrix}\frac{(-t)^j}{j!}.$$
Using the formulas for $X_+$ and $X_-$ we can see that:
$$r(X_+-X_-)\phi_2=i(2l+1)\phi_2.$$
This is because $p_l$ satisfies the differential equation:
$$p_l''(t)t+(1-t)p_l'(t)+lp_l=0.$$
Therefore the weight of the corresponding Eisenstein series is $2l+1=2(k-r)-1$. Laguerre polynomials have very nice orthogonal properties:
\begin{align} \label{orth1}
\int_0^{\infty}p_l(t)p_m(t)e^{-t}dt=\delta_{lm},
\end{align}
where $\delta_{lm}=1$ or $0$ depending on $l=m$ or not.

The other case is when the quadratic space is $(\BC,-q)$. We let $\phi_2(z)=p_l(4\pi|z|^2)e^{-2\pi|z|^2}$, where $l=r-k$ ($r\ge k$). Now the formulas for $r(X_+)$ and $r(X_-)$ on $S(\BC)$ become:
$$r(X_+)\phi=-2\pi i(z\cdot \bar{z})\phi,$$
$$r(X_-)\phi=\frac{1}{2\pi i}(\frac{\partial}{\partial
  z}\cdot\frac{\partial}{\partial {\bar{z}}})\phi.$$ 
It is easy to check that $\phi_2$ has weight $2(k-r)-1=2k-2r-1$.

\begin{remark} \label{s2,s2,r1}
In Case A2 (and Case A3 below) we have used the standard quadratic form on $\BC$, i.e, we have used $j^2=\pm1$. To get the corresponding $\phi_2$ for quadratic form $\Lambda q(z)=\Lambda |z|^2$  we only need to take $\phi_2(z)=p_l(4\pi|\Lambda||z|^2)e^{-2u\pi|\Lambda||z|^2}.$ It is not hard to see that all the statements remain true because the map $z\mapsto z/\sqrt{|\Lambda|}$ gives an isomorphism between $(\BC,q)$ and $(\BC,\Lambda q)$. 
\end{remark}

{\bf A3.}
First we let $V=\BH$, which is viewed as
$\BC\oplus\BC j$. Then we have on $S(\BH)$:
 $$r(X_+)\phi=2\pi i(x\bar{x}+y\bar{y})\phi,$$
 and
$$r(X_-)\phi=-\frac{1}{2\pi i}(\frac{\partial}{\partial
  x}\cdot\frac{\partial}{\partial {\bar{x}}}+\frac{\partial}{\partial
  y}\cdot\frac{\partial}{\partial {\bar{y}}})\phi.$$
Here $x$ and $y$ are variables on the two copies of $\BC$ respectively.
Let $\phi(z)=\phi_1(x)\phi_2(y)\in S(\BH)$, then $\phi$ and thus the corresponding  theta lift has
 weight $2r+2l+2=2k$ (which comes from a similar computation as in case A1 or A2).

Now we let $V=M_2(\BR)$, the real matrix algebra of order 2. It is also viewed as $V=\BC+\BC j$. On $S(M_2(\BR))$ we have:
$$r(X_{+})\phi=2\pi i(x\bar{x}-y\bar{y})\phi,$$
and 
$$r(X_{-})(\phi)=-\frac1{2\pi i}(\frac{\partial}{\partial x}\cdot\frac{\partial}{\partial {\bar{x}}}-\frac{\partial}{\partial
  y}\cdot\frac{\partial}{\partial {\bar{y}}})\phi.$$
Let $\phi$ still be the product of $\phi_1$ and $\phi_2$, then the corresponding theta lift has weight $2r-2l=2k$.

\begin{prop}[Archimedean weight structure] \label{s2,s2,p1}
For $\kappa=e^{i\alpha}\in SL_2(\BR)$, $k_1=e^{i\alpha_1},k_2=e^{i\alpha_2}\in K_v^{\times}\subset B^{\times}_{\BR}$ we have
$$r_B((k_1,k_2),\kappa,\phi)=e^{2ri(\alpha_1-\alpha_2)}e^{2ki\alpha}\phi.$$
\end{prop}
\begin{proof}
It is clear from the above discussions.
\end{proof}

\section{Rankin-Selberg convolution} \label{s3}

We will study the Rankin-Selberg convolution integral which involves certain Eisenstein series arising from a Weil representation. Here the Weil representation is associated to the quadratic space $(K, \Lambda\N)$ and  a fixed additive character $\psi$. Let $\delta$ be the conductor of $\psi$, i.e., $\psi(x)=\psi^0(\delta x)$ for an unramified additive character $\psi^0$ of $\BA$ (but not of $\BA/F$). Let $t(\delta)=\begin{pmatrix}\delta&0\\ 0&1\end{pmatrix}$. As before the quadratic character associated to $K/F$ is denoted by $\omega$. For $g\in G(\BA)$, $s\in\BC$, $\phi\in S(V_{\BA})$ we define a function
\begin{align} \label{s2,s3,1}
f_{}(s,g;\phi)=(r_{\Lambda}(g_1)\phi)(0)|a_{\delta}(g)|^{2s-1}|\det(g)|^{-1/2}\omega^{-1}(\det(g)),
\end{align}
where $a_{\delta}(g)=|a/b|^{1/2}$ if $g=\begin{pmatrix}a& x\\ 0& b\end{pmatrix}k'$ for $k'$  in the twisted maximal compact subgroup $U'=t(\delta)^{-1}U t(\delta)$ of $G(\BA)$, where $U$ is the standard maximal compact subgroup (as one has a twisted Iwasawa decomposition $G(\BA)=A N U'$). The function $f_{}(s,g;\phi)$ belongs to the induced representation $\mathcal{B}(|\,|^{s/2},\omega^{-1}|\,|^{-s/2})$. The related Eisenstein series is given by
$$E_{}(s,g;\phi)=\sum_{\gamma\in B(F)\cap SL_2\backslash SL_2(F)}f_{}(s,\gamma g;\phi)=\sum_{\gamma\in B(F)\backslash G(F)}f_{}(s,\gamma g;\phi).$$
This series converges for $\Re(s)\gg 0$ and  has an analytic continuation to the whole $s$-plane.

If $a(g)$ denotes the function with respect to the standard $U$, then we have the following
 \begin{equation} \label{s2,s3,2}
a_\delta(g)=a(t(\delta^{-1})gt(\delta)).
\end{equation}

\subsection{Non-archimedean case} \label{s3,s2}

Let $F=F_v$ be a nonarchimedean field. We first consider an unramified additive character $\psi^0$ of $F$. Let $\pi$  be an irreducible (infinite) representation of $G(F)=GL_2(F)$ of central character $\omega_{\pi}$. For every $W$ in the Whittaker model $\mathcal{W}(\pi,\psi^0)$ of $\pi$ we define its Mellin transform by
\begin{equation*}\Psi(s,W)=\int_{F^{\times}}W\left(\begin{pmatrix}t&0\\ 0&1\end{pmatrix}\right)|t|^{s-1/2}d^{\times}t.\end{equation*}
A Whittaker function $W_{\pi}^0\in \mathcal{W}(\pi,\psi^0)$ is called a Whittaker newform if it takes value $1$ at the unit matrix and is right invariant under $U_1(c(\pi))$, where $c(\pi)$ denotes the conductor of $\pi$.
By  results of Casselman \cite{Casselman73} and Zhang \cite{Zhang01} it is known that the Whittaker newform exists uniquely and satisfies
$$\Psi(s,W_{\pi}^0)=L(s,\pi).$$

The Whittaker newform $W_{\pi}^0$ has a  functional equation. Let $g_0=\begin{pmatrix}0&1\\-c(\pi)&0\end{pmatrix}$ and let $W_{\tilde{\pi}}^0$ be the Whittaker newform of the contragredient $\tilde{\pi}$, then (see \cite{Zhang01})
\begin{equation} \label{s3,s2,1}
W_{\pi}^0(gg_0)={W^0_{\tilde{\pi}}}(g)\omega_{\pi}(\det g)\epsilon(\pi,\psi^0),
\end{equation}
where $\epsilon(\pi,\psi^0)=\epsilon(1/2,\pi,\psi^0)$ is the local $\epsilon$-factor of $\pi$ with respect to $\psi^0$. 

Now we assume $\pi_1=\pi_v$, $\pi_2=\pi_{\chi,v}$ and let $W^0_i$ be the Whittaker newform of $\pi_i$ for $i=1$, $2$. Let $(K,\Lambda \N_K)$ be a quadratic space (splitting or not) over $F$. 

For $\phi_2\in S(K)$ we define
$$f^0(s,g;\phi_2)=r_{\Lambda}^0(g_1)\phi_2(0)|a(g)|^{2s-1}|\det(g)|^{-1/2}\omega(\det(g)),$$where  $a(g)=|a/b|^{1/2}$ if $g=\begin{pmatrix}a& x\\ 0& b\end{pmatrix}k$ with $k$ in the standard maximal subgroup $U_0(1)$. 
The modified Rankin-Selberg convolution of $W_1^0$ and $W_2^0$ given by
\begin{equation} \label{s3,s0,2}
W^{+}(s,W_1^0,W_2^0,\phi_2)=\int_{Z(F)N(F)\backslash G(F)^+}W_1^0(\epsilon g)W_2^0(g)f^0(s,g;\phi_2)dg,
\end{equation}
where $G(F)^+=\{ g\in GL_2(F):\,\det g\in \N(K)\}$.  Here the measure of $Z(F)N(F)\backslash G(F)^+$ is taken to be
\begin{equation} \label{s3,s0,3}
dg=|t|^{-1} d^{\times}t dk
\end{equation}
under the decomposition $G(F)^+=Z(F)N(F)T_1(F)U_0(1)^+$, where $d^{\times}t$ is the measure on $T_1(F)=F^{\times}$ such that $\mu(\N_K(U_K))=1$, and $dk$ on $U_0(1)^+=U_0(1)\cap G^+$ is such that $\mu(U_0(1)^+)=1$.

\begin{prop} \label{s3,s2,p1}
Let $c={\rm{Max}}(v(c(\pi_1)),v(c(\pi_2)))$. Let $W_1^0\in\mathcal{W}(\pi_1,\psi^0)$ and $W_2^0\in \mathcal{W}(\pi_2,\psi^0)$ be the Whittaker newforms of $\pi_1$ and $\pi_2$  respectively, then

\begin{equation*}
L(s,\pi_1\times\pi_2)=\begin{cases}
{L(2s,\omega)}\Psi^+(s,W^0_1,W^0_2,\phi_2)\quad & \text{in cases (1) and (2),}\\
\frac{L(1,\omega)}{\mu(U_0(\varpi^c)^+)}\Psi^+(s,W^0_1,W^0_2,\phi_2)\quad & \text{in cases (3), (4), (5).}
\end{cases}
\end{equation*}
where $\omega$ as usual denotes the quadratic character of $F$ associated to $K/F$, and the five cases are divided according to Lemma \ref{s2,s3,l1}.
\end{prop}
\begin{proof}
Most of the claims are already proved in \cite{Popa03} (for squarefree $N$). We will supply a simpler proof for some cases. The method, which  is first used in \cite{Zhang01},  can also  treat non-squarefree $N$. 

Cases 1, 2 and 5 can be proved using Corollary \ref{s3,s2,c1} Cases 1,2 and 5 respectively. Since the proof is the same as the one used in Cases A1, A1' and A3 of \cite{Popa03}, so is omitted.

Case 3a occurs when $c=v(c(\pi_1))=2t+1$ and $v(\Lambda)=1$. As $W_1^0$ has the minimum level $U_1(\varpi^c)$ and $\pi_1$ has conductor $\varpi^c$, so 
$$\sum_{\gamma\in U_1(\varpi^{c-1})/U_1(\varpi^c)}W_1^0(g\gamma)=0.$$
Then 
\begin{equation} \label{s3,s2,-2}
\int W_1^0(\epsilon g) W(g)dg=\frac{1}{|U_1(\varpi^{c-1})/U_1(\varpi^c)|}\sum_{\gamma\in U_1(\varpi^{c-1})/U_1(\varpi^c)}\int W_1^0(\epsilon g)W(g\gamma^{-1})dg=0
\end{equation}
for any function $W$ invariant under $U_1(\varpi^{c-1})$.

By Lemma \ref{s3,s2,l2.5} $f^0(s,k;\phi_2)$ is the sum of $(1+|\varpi|)1_{U_0(\varpi^c)}$ and a function which is invariant under $U_0(\varpi^{c-1})$, so
\begin{align*}
\Psi^+(s,W_1^0,W_2^0,\phi_2)=&(1+|\varpi|)\int_{T_1^+(F)\times U_0(1)}W_1^0(-tk)W_2^0(tk)1_{U_0(\varpi^c)}|t|^{s-1}dk dt\\
=&\frac{\mu(U_0(\varpi^c))}{L(1,\omega)}\int_{T_1^(F)}W_1^0(-t)W_2^0(t)|t|^{s-1}dt\\
=&\frac{\mu(U_0(\varpi^c))}{L(1,\omega)}L(s,\pi_1\times\pi_2),
\end{align*}
where the last equality is a well-known fact. Since $U_0(\varpi^c)=U_0(\varpi^c)^+$ the proof of case 3a is complete. 

Case 3b occurs when either $v(c(\pi_1))=c=2t$ or $v(c(\pi_2))=c=2t$. This case can be treated in the same way as case 3a.

Case 4 shows up when $v(c(\pi_2))=c$. As $L(1,\omega)^{-1}=(1-|\varpi|)$ the above argument for Case 3a moves parallelly to this case (also using Lemma \ref{s3,s2,l2.5} for Case 4).
\end{proof}

\begin{lem} \label{s3,s2,l1}
If $K/F$ is a quadratic field extension of different $\pi_K^C\mathcal{O}_K$, and $\alpha\in F^{\times}$, then
\begin{equation*}
\int_{\pi_K^t\mathcal{O}_K}\psi^0(\alpha y\bar{y})dy=\begin{cases}|\pi|^{ft}\mu(\mathcal{O}_K) & \text{if $v(\alpha)\ge-ft$,}\\
0 & \text{if $-ft>v(\alpha)>-C-ft$,}\\
\gamma|\alpha|^{-1}\omega(\alpha)^{-1} & \text{if $v(\alpha)\le-C-ft$,}
\end{cases}
\end{equation*}
here $f$ is the residue degree of $K/F$.

\noindent If $K=F\oplus F$ is split, then 
\begin{equation*}
\int_{\pi^s\mathcal{O}\times\pi^t\mathcal{O}}\psi^0(\alpha xy)dxdy=\begin{cases}|\pi|^{s+t} & \text{if $v(\alpha)\ge-s-t$,}\\
|\alpha|^{-1} & \text{if $v(\alpha)\le-s-t$.}
\end{cases}
\end{equation*}
\end{lem}
\begin{proof}
This is Lemma 2.5.2 in \cite{Popa03}.
\end{proof}

\begin{lem} \label{s3,s2,l2}
Let $k=\begin{pmatrix}a&b\\ c&d\end{pmatrix}\in U_0(1)$ and $v(\Lambda)=M$.
If $K/F$ is inert  and $\phi_2=1_{\varpi_K^{\ell}\mathcal{O}_K}$, then
\begin{equation*}
f^0(s,k;\phi_2)=\begin{cases}1\quad & \text{if $k\in U_0(\varpi^{M+2\ell})$,}\\
|\Lambda c^{-1}\varpi^{2\ell}|\omega(\Lambda c^{-1})\quad & \text{otherwise.}
\end{cases}
\end{equation*}

If $K=F\oplus F$ and $\phi_2=1_{\mathcal{O}\times\varpi^{\ell}\mathcal{O}}$, then
\begin{equation*}
f^0(s,k;\phi_2)=\begin{cases}1\quad & \text{if $k\in U_0(\varpi^{M+\ell})$,}\\
|\Lambda c^{-1}\varpi^{\ell}|\quad & \text{otherwise.}
\end{cases}
\end{equation*}

If $K/F$ is ramified of ramification index $1$ and $\phi_2=1_{\varpi_K^{-M}\mathcal{O}_K}$, then
\begin{equation*}
f^0(s,k;\phi_2)=\begin{cases}\omega(d)\quad & \text{if $k\in U_0(\varpi)$,}\\
\gamma\mu(\mathcal{O}_K)|\varpi^{-M}\Lambda c^{-1}|\omega(-\Lambda c^{-1})\quad & \text{otherwise.}
\end{cases}
\end{equation*}
\end{lem}
\begin{proof} The claims can be proved through  direct computations using explicit formulae in Section 2.1, Lemma  \ref{s3,s2,l1}  and  the Bruhat decomposition 
$$G=B\coprod B\begin{pmatrix}0&1\\-1&0\end{pmatrix}N.$$ More precisely, if $c\ne0$ and $a\ne0$ we have
\begin{align*}
\begin{pmatrix}a&b\\c&d\end{pmatrix}=-\begin{pmatrix}c^{-1}&a\\0&c\end{pmatrix}\begin{pmatrix}0&1\\-1&0\end{pmatrix}\begin{pmatrix}1&ba^{-1}+c^{-1}a^{-1}\\0&1\end{pmatrix},
\end{align*}
while if $c\ne0$ and $a=0$ we have
\begin{align*}
\begin{pmatrix}0&b\\c&d\end{pmatrix}=\begin{pmatrix}0&1\\-1&0\end{pmatrix}\begin{pmatrix}-c&-d\\0&b\end{pmatrix}.
\end{align*}
The detail is omitted.
\end{proof}

\begin{cor} \label{s3,s2,c1}
According to the five cases of Proposition \ref{s3,s2,p1}:\\
(1):  $\phi_2=1_{\varpi^{-m}\mathcal{O}_K}$, then
$f^0(s,k;\phi_2)=1$;\\ 
(2):  $\phi_2=1_{\mathcal{O}}\times 1_{\varpi^{-m}\mathcal{O}}$, then $f^0(s,k;\phi_2)=1$;\\
(3):  $\phi_2=1_{\varpi^t\mathcal{O}_K}$, then 
\begin{equation*}
f^0(s,k;\phi_2)=\begin{cases}1\quad & \text{if $k\in U_0(\varpi^c)$,}\\
|\Lambda c^{-1}|\omega(\Lambda c^{-1})\quad & \text{otherwise;}
\end{cases}\end{equation*}
(4):  $\phi_2=1_{\mathcal{O}}\times1_{\varpi^{c-m}\mathcal{O}}$ for $m=0,1$, then  
\begin{equation*}
f^0(s,k;\phi_2)=\begin{cases}1\quad & \text{if $k\in U_0(\varpi^c)$,}\\
|\Lambda c^{-1}| \quad &\text{otherwise;}
\end{cases}
\end{equation*}
(5):  $\phi_2=1_{\mathcal{O}_K}$, then
\begin{equation*}
f^0(s,k,\phi_2)=\begin{cases}\omega(d)\quad & \text{if $k\in U_0(\varpi)$,}\\
\gamma\mu(\mathcal{O}_K)|\Lambda|\omega(-\Lambda c^{-1}) \quad &\text{otherwise.}\end{cases}
\end{equation*}
here the norm $|\cdot|$ is the norm in $F$, and in  case (5) we have assumed the residue characteristic of $F$ is odd.

Therefore,
$f^0(s,g;\phi_2)$ is right invariant under $U_0(\varpi^{c})$ in the first 4 cases and  is right invariant under $U_0(\varpi^{c})^+$ in Case 5.
\end{cor}

\begin{lem} \label{s3,s2,l2.5}
Let $k\in U_0(1)$.
In Case 3,  one has $$f^0(s,k;\varphi_2)=(1+|\varpi|)1_{U_0(\varpi^{2t})}-|\varpi|(1+|\varpi|)1_{U_0(\varpi^{2t-1})}+f(s,k;|\varpi|^21_{\varpi^{t-1}\mathcal{O}_K}).$$

\noindent In Case 4, one has $$f^0(s,k;\phi_2)=(1-|\varpi|)1_{U_0(\varpi^{c})}+|\varpi|f^0(s,k;1_{\mathcal{O}\times\varpi^{c-m-1}\mathcal{O}}).$$
\end{lem}
\begin{proof}
For Case 3 by Lemma \ref{s3,s2,l2}  
$$f^0(s,k,1_{\varpi^{t-1}\mathcal{O}_K}-|\varpi|^{-2}1_{\varpi^{t}\mathcal{O}_K})=(|\varpi|\omega(\varpi)-|\varpi|^{2})1_{U_0(\varpi^{2t-1})}-(1-|\varpi|\omega(\varpi))1_{U_0(\varpi^{2t})},$$
hence the conclusion.

For Case 4 Lemma \ref{s3,s2,l2} again implies 
$$f^0(s,k;\phi_2-|\varpi| 1_{\mathcal{O}\times\varpi^{c-m-1}})=(1-|\varpi|)1_{U_0(\varpi^c)},$$
from which the result follows.
\end{proof}

When $v$  divides $c(\chi)^2=c(\pi_{2})$, or equivalently in Cases 3b and 4, we need a stronger result similar to  Proposition 2.5.1 in \cite{Zhang01}.

\begin{prop} \label{s3,s2,p2}

Assume $0\le j\le 2t-1=2c(\pi_2)-1$, then
$$\Psi^+\left(s,\rho\begin{pmatrix}\varpi^{-j}&0\\0&1\end{pmatrix}W^0_1,W^0_2;\phi_2\right)=|\varpi|^{j(s-1/2)}\alpha_j\frac{\mu(U_0(\varpi^c)}{L(1,\omega)}L(s,\pi_1\times\pi_2),$$
where $\alpha_j$ is defined by
$$L(s,\pi_2)=\sum_{j=0}^{\infty}\alpha_j|\varpi|^{js}=1.$$
\end{prop} 
\begin{proof}
As $\pi_2$ has conductor $\varpi^{2t}$, 
$$\sum_{\gamma\in U_1(\varpi^{2t-1})/U_1(\varpi^{2t})}\rho(\gamma)W_2^0=0.$$

If $j<2t$, then $\rho(t(\varpi^{-j}))W_1^0$ is invariant under $U_1(\varpi^{2t-1})$. By Lemma \ref{s3,s2,l2.5} $f^0(s,k;\phi_2)$ is the sum of $L(1,\omega)^{-1}1_{U_0(\varpi^c)}$ and a function which is invariant under $U_0(\varpi^{c-1})$. So, by (\ref{s3,s2,-2}) one has  
\begin{align*}
\Psi^+(s,\rho(t(\varpi^{-j})W_1^0,W_2^0,\phi_2)&=\frac{\mu(U_0(\varpi^c))}{L(1,\omega)}\int_{T_1^+}W_1^0(t(-a\varpi^{-j}))W_2^0(t(a))|a|^{s-1}da^{\times}\\
&=|\varpi|^{(s-1/2)j}\frac{\mu(U_0(\varpi^c))}{L(1,\omega)}\alpha_j\int_{T_1^+}W_1^0(t(-a))W_2^0(t(a\varpi^j))|a|^{s-1}da^{\times}\\
&=|\varpi|^{j(s-1/2)}\alpha_j\frac{\mu(U_0(\varpi^c))}{L(1,\omega)}L(s,\pi_1\times\pi_2),
\end{align*}
where we have used the fact that $W_2^0(t(a\varpi^j))=|\varpi|^{j/2}\alpha_jW_2^0(t(a))$ for any integral $a$.
\end{proof}

We assume now $j=2t=c$. This time it is hard to prove a similar result above, but one can still obtain the information on the convolution at $s=1/2$, which is sufficient for later applications. First, the functional equations (\ref{s3,s2,1}) for $W_1^0$ and $W_2^0$ imply that (as $\pi_1=\tilde{\pi}_1$, $\omega_{\pi_1}=1$ and $\omega$ is unramified)
$$W_1^0(gt(\varpi^{-2t}))=\epsilon(\pi_1,\psi^0)W_1^0\left(gt(\varpi^{-2t})\begin{pmatrix}0&1\\-1&0\end{pmatrix}\right)=\epsilon(\pi_1,\psi^0)W_1^0\left(g\begin{pmatrix}0&1\\-\varpi^{2t}&0\end{pmatrix}\right),$$
$$W_2^0\left(g\begin{pmatrix}0&1\\-\varpi^{2t}&0\end{pmatrix}\right)=\epsilon(\pi_2,\psi^0)\tilde{W}_2^0(g)\omega(\det g)=\epsilon(\pi_2,\psi^0)W_2^0(g),$$
where $\tilde{W}_2^0$ is the Whittaker newform for $\tilde{\pi}_2$. Write $g_0=\begin{pmatrix}0&1\\-\varpi^{2t}&0\end{pmatrix}$ and $C=\epsilon(\pi_1,\psi^0)\epsilon(\pi_2,\psi^0)$, then 
\begin{align} \label{s3,s2,2.5}
&\Psi^+(s,\rho(t(\varpi^{-2t})W_1^0,W_2^0,\phi_2)\\=&\epsilon(\pi_1,\psi^0)\Psi^+\left(s,\rho(g_0)W_1^0,{W}_2^0,\phi_2\right)\nonumber\\
=&C\int_{Z(F)N(F)\backslash G(F)^+} W_1^0(\epsilon g){W}_2^0(g)f^0(s,gg_0;\phi_2)dg.\nonumber
\end{align}

\begin{lem} \label{s3,s2,l3.5}
Let $h\in K$ such that $\N(h)=\det g_0=\varpi^{2t}$ and define
\begin{equation} \label{s3,s2,-1}
\phi_2'=r^0_{\Lambda}(h,g_0)\phi_2,
\end{equation}
then 
\begin{align} \label{s3,s2,-3}
f^0(1/2,g;\phi_2')=f^0\left(1/2,gg_0;\phi_2\right),
\end{align}
and
$f^0(s,g;\phi_2')$
 is right invariant under $U_0(1)$.
\end{lem}
\begin{proof}
The identity (\ref{s3,s2,-3}) is easily verified  through a simple calculation. 
To ease the computation  we take $h=\varpi^{t}$ in Case 3b and take $h=(\varpi^{2t},1)$ in Case 4.   By definition
$$r_{\Lambda}^0\left(h,g_0\right)\phi_2(x)=L(h)r_{\Lambda}^0(w)\phi_2(x).$$
In Case 3b, $\phi_2=1_{\varpi^t\mathcal{O}_K}$,  since $v(\Lambda)=0$ and $K/F$ is unramified, we get
$$r_{\Lambda}^0(w)\phi_2=1_{\varpi^{-t}\mathcal{O}}.$$
Hence $r_{\Lambda}^0\left(h,g_0\right)\phi_2=L(h)1_{\varpi^{-t}\mathcal{O}}=|\varpi^t|1_{\mathcal{O}}.$ 

In Case 4, $K=F\oplus F$, $\phi_2=1_{\mathcal{O}\times\varpi^{2t-m}\mathcal{O}}$ and $v(\Lambda)=m$, so
$$r_{\Lambda}^0(w)\phi_2=|\Lambda|1_{\varpi^{-2t}\mathcal{O}\times\varpi^{-m}\mathcal{O}}.$$
Therefore, $r_{\Lambda}^0\left(h,g_0\right)\phi_2=L(h)1_{\varpi^{-2t}\mathcal{O}\times\varpi^{-m}\mathcal{O}}=|\Lambda|1_{\mathcal{O}\times\varpi^{-m}\mathcal{O}}.$ 

The claim now follows from Corollary \ref{s3,s2,c1} and Lemma \ref{s3,s2,-2}. 
\end{proof}
\begin{prop} \label{s3,s2,p3} 
We have
\begin{equation*}
\Psi^+(s,W_1^0,W_2^0,\phi'_2)=0,
\end{equation*}
\end{prop}
\begin{proof}
As both $W_1^0$ and $f^0(s,g;\phi')$ are invariant under $U_0(1)$, thus the claim is clear by (\ref{s3,s2,-2}) and (\ref{s3,s2,2.5}).
\end{proof}

We now assume  that $\psi$ is a general additive character  given by $\psi(x)=\psi^0(\delta x)$. 
 Let $r_{\Lambda}$ be the Weil representation associated to $(K,\Lambda\N)$ with respect to $\psi$. It is related to $r^0_{\Lambda}$ for $g_1\in SL_2(F)$ through 
\begin{equation} \label{s3,s2,1} r_{\Lambda}(g_1)=r_{\Lambda}^0\left(t(\delta)g_1t(\delta)^{-1}\right).
\end{equation}

 Let  $W_1(g)=W^0_1(t(\delta)gt(\delta^{-1}))$ and $W_2(g)=W^0_2(t(\delta)gt(\delta^{-1}))$, so they are in  $\mathcal{W}(\pi_1,\psi)$ and $\mathcal{W}(\pi_2,\psi)$ respectively. The modified Rankin-Selberg convolution between $W_1$ and $W_2$ is given by
\begin{equation}
W^{+}(s,W_1,W_2,\phi_2)=\int_{Z(F)N(F)\backslash G^{+}(F)}W_1(\epsilon g)W_2(g)f(s,g;\phi_2)dg,
\end{equation} 
where $f(s,g;\phi_2)=r_{\Lambda}(g_1)\phi_2(0)|a_{\delta}(g)|^{2s-1}|\det(g)|^{-1/2}\omega(\det(g))$. Here, the measure $dg$  on $Z(F)N(F)\backslash G^{+}(F)$ is induced by that of (\ref{s3,s0,3}) through the automorphism $g\to t(\delta^{-1})gt(\delta)$, and it is easy to see they are the same. 

By (\ref{s2,s3,2}) and (\ref{s3,s2,1}) 
\begin{equation} \label{s3,s1,1}
f(s,t(\delta)^{-1}gt(\delta);\phi_2)=f^0(s,g;\phi_2).
\end{equation}

\begin{prop} \label{s3,s2,p4}
Let $W_i(g)=W^0_i(t(\delta)gt(\delta^{-1}))$ for $i=1$, 2. For $\Re(s)\gg0$
\begin{equation} \label{s3,s2,3}
\Psi^+(s,W_1,W_2,\phi_2)=\Psi^+(s,W^0_1,W^0_2,\phi_2).
\end{equation}
For $1\le j< 2t=2v(c(\chi))$ 
\begin{equation} \label{s3,s2,4}
\Psi^+\left(s,\rho(t(\varpi^{-j}))W_1,W_2,\phi_2\right)=0.
\end{equation}
\end{prop}
\begin{proof}
Let $g'=t(\delta) g t(\delta^{-1})$, then 
\begin{align*}
&\Psi^+(s,W_1,W_2,\phi_2)\\
=&\int_{Z(F)N(F)\backslash G^+(F)}W_1(\epsilon t(\delta^{-1})g't(\delta))W_2(t(\delta^{-1})g't(\delta))f(s,t(\delta^{-1})g't(\delta);\phi_2)dg'\\
=&\int_{Z(F)N(F)\backslash G^+(F)}W^0_1(\epsilon g')W^0_2(g')f^0(s,g';\phi_2)dg',
\end{align*}
which is exactly the right hand side of (\ref{s3,s2,3}).

The proof of (\ref{s3,s2,4})  is similar and uses Proposition \ref{s3,s2,p2}.
\end{proof}

\begin{prop} \label{s3,s2,p5}
For $j=2v(c(\chi))>0$ and $C=\epsilon(\pi_1,\psi^0)\epsilon(\pi_2,\psi^0)$
\begin{equation} \label{s3,s2,5}
\Psi^+\left(s,\rho(t(\varpi^{-j}))W_1,W_2,\phi_2\right)=C\int_{Z(F)N(F)\backslash G(F)^+} W_1(\epsilon g){W}_2(g)f(s,gt(\delta^{-1})g_0t(\delta);\phi_2)dg,
\end{equation}
and for $\phi_2'$ in (\ref{s3,s2,-1})
\begin{equation} \label{s3,s2,6}
f(1/2,gt(\delta^{-1})g_0t(\delta);\phi_2)=f(1/2,g;\phi_2').
\end{equation}
Moreover,
\begin{equation} \label{s3,s2,7}
\Psi^+(s,W_1,W_2,\phi'_2)=0.
\end{equation}
\end{prop}
\begin{proof}
Let $g'=t(\delta)gt(\delta^{-1})$. The left side of (\ref{s3,s2,5}) becomes
\begin{align*}
&\int_{Z(F)N(F)\backslash G(F)^+}W_1(\epsilon gt(\varpi^{-j}))W_2(g)f(s,g;\phi_2)dg\\
=&\int_{Z(F)N(F)\backslash G(F)^+}W^0_1(\epsilon g't(\varpi^{-j}))W^0_2(g')f^0(s,g';\phi_2)dg'\\
=&C\int_{Z(F)N(F)\backslash G(F)^+} W_1^0(\epsilon g'){W}_2^0(g')f^0(s,g'g_0;\phi_2)dg',
\end{align*}
where in the last step we have used (\ref{s3,s2,2.5}).
The right side of (\ref{s3,s2,5}) becomes
\begin{align*}
&C\int_{Z(F)N(F)\backslash G(F)^+} W_1(\epsilon g){W}_2(g)f(s,gt(\delta^{-1})g_0t(\delta);\phi_2)dg\\
=&C\int_{Z(F)N(F)\backslash G(F)^+} W^0_1(\epsilon g'){W}^0_2(g')f(s,t(\delta^{-1})g'g_0t(\delta);\phi_2)dg'\\
=&C\int_{Z(F)N(F)\backslash G(F)^+} W^0_1(\epsilon g'){W}^0_2(g')f^0(s,g'g_0;\phi_2)dg'.
\end{align*}
By (\ref{s3,s2,-3})
\begin{align*}
f(1/2,gt(\delta^{-1})g_0t(\delta);\phi_2)&=f^0(1/2,t(\delta)gt(\delta^{-1})g_0t(\delta)t(\delta^{-1});\phi_2)\\
&=f^0(1/2,t(\delta)gt(\delta^{-1});\phi'_2).
\end{align*}
The proof of (\ref{s3,s2,7}) is the same as that of (\ref{s3,s2,3}).
\end{proof}

Over a finite place $v|c(\chi)$ we now pick a special vector in $\mathcal{W}(\pi_1,\psi)$. The  space $\mathcal{W}(\pi_1,\psi)$ is equipped   with an inner product, denoted  by $(\, ,\,)$. Let $W_1 \in\mathcal{W}(\pi_1,\psi)$ be the Whittaker function defined in Proposition \ref{s3,s2,p4}. By  newform theory (see Lemma \ref{s3,s2,l4} below for a proof in our twisted situation), the subspace of vectors in $\mathcal{W}(\pi_1,\psi)$ which are invariant under $t(\delta^{-1})U_1(c(\chi)^2) t(\delta)$ is generated by $\rho(t(\varpi^{-j}))W_1$ for $j=0,\cdots,2v(\chi).$ 

We define  $W_1^{\ast}$ to be a vector  invariant under $t(\delta^{-1})U_1(c(\chi)^2) t(\delta)$, such that

{\bf 1.} $(W_1^{\ast},W_1-W_1^{\ast})=0$,

{\bf 2.} $\left(W_1^{\ast},\rho\begin{pmatrix}\varpi^{-j}&0\\0&1\end{pmatrix}W_1\right)=0$ for $0<j\le v(c)$. 

\noindent Such a vector is called a local quasi-newform in \cite{Zhang01}.
Note that $W_1^{\ast}$ is unique and does not depend on the choice of inner product. 

\begin{lem} \label{s3,s2,l4}
We have $W_1^{\ast}(1)=1$.
\end{lem}
\begin{proof}
 Let $\tilde{W}_1^{\ast}=W_1^{\ast}(t(\delta)^{-1}gt(\delta))$, so $\tilde{W}_1^{\ast}\in\mathcal{W}(\pi_1,\psi^0)$ and is invariant under $U_1(c(\chi)^2)$. By  newform theory (\cite{Casselman73} \cite{Zhang01}) $\tilde{W}^{\ast}_1$ is of the form
 $$\tilde{W}^{\ast}_1=\sum_{j=0}^{2t}c_j\rho(t(\varpi^{-j}))W^0_1.$$
Conjugating back by $t(\delta)$ 
$$W_1^{\ast}=\sum_{j=0}^{2t}c_j\rho(t(\varpi^{-j}))W_1.$$
 So 
$$W^{\ast}_1(1)=\sum_{j=0}^{2t}c_jW_1^0(t(\varpi^{-j}))=c_0,$$ 
and 
$$(W^{\ast}_1,W_1^{\ast})=(W^{\ast}_1,c_0W_1+\sum_{j\ne0}c_j\rho(t(\varpi^{-j}))W_1))=\bar{c}_0(W^{\ast}_1,W_1),$$ therefore $W_1^{\ast}(1)=c_0=1$. 
\end{proof}

\subsection{Archimedean case} \label{s3,s1}
Assume now that $F=\BR$  and let $\psi$ be the
 standard  additive character of $F$ defined by $\D \psi(x)=e^{2\pi ix}$. We write $\pi_1=\pi_{v}$ and $\pi_2=\pi_{\chi,v}$.

The representation $\pi_1$ is a discrete series of weight $2k=2k_v$ and $\pi_2$ is a discrete series of weight $2r+1=2r_v+1$. Let $\pi$ be a discrete series of lowest weight $k$, then  $W^0\in\mathcal{W}(\pi,\psi)$ is called a standard Whittaker function of weight  $k$ if its weight is $k$, $W^0(t(a))=0$ for $a<0$ and
$$\int_{\BR}W^0(t(a))|a|^{s-1/2}d^{\times}a=L(s,\pi).$$
Such a Whittaker function satisfies
\begin{equation*}
W^0(t(a))=\begin{cases}2a^{k/2}e^{-2\pi a} & \text{if $a>0$,}\\
0 & \text{if $a<0$.}\end{cases}
\end{equation*}
 Similarly, we say $W^0$ is a standard Whittaker function of weight $-k$ if  its weight is $-k$, and
\begin{equation*}
W^0(t(a))=\begin{cases}2|a|^{k/2}e^{2\pi a} & \text{if $a<0$,}\\
0 & \text{if $a>0$.}\end{cases}
\end{equation*}

The definition of $L(s,\pi)$ will be recalled during the proof of the following proposition.

 Locally at $v$ (we drop $v$ occasionally in the rest of this section):
$$f_{}(s,g;\phi_2)=r_{\Lambda}(g_1)\phi_2(0)|a(g)|^{2s-1}|\det(g)|^{-1/2}\omega(\det(g)),$$
where $\phi_2$ is defined in  Remark \ref{s2,s2,r1}.
Note that $\phi_2(0)=1$, so $f(s,1,\phi_2)=1$. In the
following we take $l=k-r-1$ if $k>r$  and $l=r-k$ if $r\ge k$. As in Section \ref{s1} we write $\Sigma_1$ for the set of  infinite $v$ with $k_v>r_v$ and $\Sigma_2$ for the set of $v$ with $r_v\ge k_v>0$.

\begin{lem}[Barnes' Lemma] \label{s3,s1,l1}
Assume $f_1$, $f_2$ are smooth functions on $\BR_+$,  such that
$$\int_{\BR_+}f_1(a)|a|^{s-1/2}d^{\times}a=G_1(s+r_1)G_1(s+r_2)$$
$$\int_{\BR_+}f_2(a)|a|^{s-1/2}d^{\times}a=G_1(s+t_1)G_1(s+t_2),$$
then
$$\int_{\BR_+}f_1(a)f_2(a)|a|^{s-1}d^{\times}a=\frac{\prod_{i,j=1}^2G_1(s+r_i+t_j)}{G_1(2s+r_1+r_2+t_1+t_2)}.$$
\end{lem}

\begin{prop} \label{s3,s1,p1}
Let $W_1$  be the standard Whittaker function of weight $-2k$  in $\pi_1$ and let $W_2$ be the standard Whittaker function of weight $2r$ in $\pi_2$.  Then
\begin{align*}
L(s,\pi_1\times\pi_2)=\begin{cases}2^{s+k+r+1/2}G_2(s+k-r-1/2)\cdot\Psi^+(s,W_1,W_2,\phi_2),\quad &\text{if $v\in\Sigma_1$,}\\
2^{s+k+r+1/2}G_2(s+r-k+1/2)\cdot\Psi^+(s,W_1,W_2,\phi_2),\quad &\text{if $v\in\Sigma_2$.}
\end{cases}
\end{align*}
Here the measure is taken such that the total measure of $SO_2(\BR)$ is $1$.
\end{prop}
\begin{proof}
By the Iwasawa decomposition and invariance of the triple product under $SO_2(\BR)$ from the weight computation 
\begin{align*}
\Psi^+(s)=&\Psi^+(s,W_1,W_2,\phi_2)=\int_{Z(\BR)N(\BR)\backslash G(\BR)^+}W_1(\epsilon g)W_2(g)f(s,g;{\phi_2})dg\\
=&\int_{\BR^{+}}W_1\left(\begin{pmatrix}a&0\\ 0&1\end{pmatrix}\right)W_2\left(\begin{pmatrix}a&0\\ 0&1\end{pmatrix}\right)|a|^{s-1}f(s,1,\phi_2)d^{\times}a.
\end{align*}

Recall that
$$G_1(s)=\pi^{-s/2}\Gamma(s/2),\quad G_2(s)=2(2\pi)^{-s}\Gamma(s).$$

If $v$ is in $\Sigma_1$, i.e. $k>r$, then the local $L$-factors are given by
$$L(s,\pi_1)=G_2(s+k-\frac12),\quad L(s,\pi_2)=G_2(s+r),$$
and
\begin{align*}
L(s, \pi_1\times\pi_2)=G_2(s+k+r-1/2)\cdot G_2(s+k-r-1/2)
\end{align*}
Now we  let $r_1=k-1/2$, $r_2=k+1/2$, $t_1=r$, $t_2=r+1$ and apply the Barnes' formula:
$$\Psi(s)=\frac{G_1(s+k+r-1/2)G_1(s+k+r+1/2)G_1(s+k+r+1/2)G_1(s+k+r+3/2)}{G_1(2s+2k+2r+1)}.$$
So the difference between $\Psi$ and the Rankin-Selberg convolution is
given by
\begin{align*}
L(s,\pi_1\times\pi_2)&=\Psi(s)\cdot
\frac{G_2(s+k-r-1/2)G_1(2s+2k+2r+1)}{G_2(s+k+r+1/2)}\\
&=2^{s+k+r+1/2}G_2(s+k-r-1/2)\cdot \Psi(s),
\end{align*}
where we have used the formula
$$G_1(2s)=2^{s-1}G_2(s)=2^{s-1}G_1(s)G_1(s+1).$$
The  case when $v\in\Sigma_2$ can be treated exactly in the same way.
\end{proof}

\subsection{Global case} \label{s3,s3}
 Let $\pi_1=\pi=\prod_v\pi_{v}$ and $\pi_2=\pi_{\chi}=\prod_v\pi_{\chi,v}$. 
We take the measure on $G(\BA)^+$  to be the product of the local ones used in previous sections, and define
$$W=\prod_{v}W_{1,v},\,\,\, W^{\ast}=\prod_{v\not|c(\chi)}W_{1_v}\prod_{v\in c(\chi)}W_{1,v}^{\ast},\,\,\, W_{\chi}=\prod_v W_{2,v},$$ 
where these local Whittaker functions are given in Sections \ref{s3,s2} and \ref{s3,s1}. 
The corresponding automorphic forms on $G(\BA)$ are constructed  as follows
\begin{align} \label{s3,s3,0}
\varphi(g)=\sum_{\xi\in F^{\times}}W\left(\begin{pmatrix}\xi&0\\0&1\end{pmatrix}g\right),
\end{align}
\begin{align} \label{s3,s3,1} 
\varphi^{\ast}(g)=\sum_{\xi\in F^{\times}}W^{\ast}\left(\begin{pmatrix}\xi&0\\0&1\end{pmatrix}g\right),
\end{align}
\begin{align} \label{s3,s3,2}
\varphi_{\chi}(g)=C_{\chi}(g)+\sum_{\xi\in F^{\times}}W_{\chi}\left(\begin{pmatrix}\xi&0\\0&1\end{pmatrix}g\right),
\end{align}
where $C_{\chi}(g)$ is the (possibly zero) constant term of $\varphi_{\chi}$. 
For a nonempty set $T$ of places dividing $c(\chi)$ we let $g^T_0=\prod_{v\in T}g_{0,v}$, where $g_{0,v}=\begin{pmatrix}0&1\\-c(\chi_v)^2&0\end{pmatrix}$. We also write 
$$\phi_2=\prod_v\phi_{2,v},\quad \phi_{2,T}'=\prod_{v\not\in T}\phi_{2,v}\prod_{v\in T}\phi'_{2,v},$$ 
where $\phi_{2,v}\in S(K_v)$ is given in Propositions \ref{s3,s2,p1} and \ref{s3,s1,p1}, and $\phi_{2,v}'$ is defined in Lemma \ref{s3,s2,l3.5}. 
The Eisenstein series associated to $\phi_2$  is defined by
\begin{align}
E(s,g;\phi_2)=\sum_{\gamma\in B(F)\backslash G(F)}f(s,\gamma g;\phi_2),
\end{align}
and the Eisenstein series attached to $\phi'_2$ is given  by
\begin{align}
E(s,g;\phi'_{2,T})=\sum_{\gamma\in B(F)\backslash G(F)}f(s,\gamma g;\phi'_{2,T}),
\end{align}
where $f(s,g;\phi)=r_{\Lambda}(g_1)\phi(0)|a_{\delta}(g)|^{2s-1}|\det g|^{-1/2}\omega(\det g)$ for any $\phi\in S(K_{\BA})$. Both Eisenstein series have analytic continuation over the whole $s$-plane.  By (\ref{s3,s2,6}) and analytic continuation
\begin{equation} \label{s3,s3,2.5}
E(1/2,g t(\delta^{-1})g_0^Tt(\delta);\phi_2)=E(1/2,g;\phi'_{2,T}).
\end{equation}

\begin{lem} \label{s3,s3,l1}
Let $a|c(\chi)^2=c(\pi_{\chi})$ be a  non-unit integral idele and let  $\varphi_a=\rho(t(a^{-1}))\varphi$, then
\begin{equation}
\int_{Z(\BA)G(F)^+\backslash G(\BA)^+}\varphi_a(g)\varphi_{\chi}(g)E(1/2,g;\phi_2)dg=0.
\end{equation}
\end{lem}
\begin{proof}
We write $W_a=\rho(t(a^{-1}))W$. If $v(a)<c(\pi_{\chi})$ for all $v$, then  
\begin{align*}
&\int_{Z(\BA)G(F)^+\backslash G(\BA)^+}\varphi_a(g)\varphi_{\chi}(g)E(s,g;\phi_2)dg\\=&\int_{Z(\BA)N(\BA)\backslash G(\BA)^+}W_a(\epsilon g)W_{\chi}(g)f(s,g;\phi_2)dg=0.
\end{align*}
If $v(a)=v(c(\pi_{\chi}))$ for $v\in T$,  then by Proposition \ref{s3,s2,p4}
\begin{align*}
 \int_{Z(\BA)G(F)^+\backslash G(\BA)^+}\varphi_a(g)\varphi_{\chi}(g)E(1/2,g;\phi_2)dg
=\int_{Z(\BA)N(\BA)\backslash G(\BA)^+}\varphi_a(g)\varphi_{\chi}(g)E(1/2,g;\phi'_{2,T})dg,
\end{align*}
which is zero by  (\ref{s3,s2,7}) through analytic continuation.
\end{proof}

We say $\psi\in\pi$ is holomorphic  (or anti-holomorphic) if its Whittaker function 
$$\D W_{\psi}(g)=\int_{F\backslash \BA}\psi\left(\begin{pmatrix}1&x\\0&1\end{pmatrix}g\right)dx$$ satisfies
$$W_{\psi}\left(\begin{pmatrix}a_{\infty}a_f&0\\0&1\end{pmatrix}\right)=\widehat{\psi}(a_f)W_{\infty}(t(a_{\infty}))\,\,\,( {\rm or}\,\,\widehat{\psi}(a_f)W_{\overline{\infty}}(t(a_{\infty}))),$$
where $\widehat{\psi}$ is a function of finite ideles $a_f$, $W_{\infty}=\prod_{v|\infty}W_v$ (or $W_{\overline{\infty}}$) is the standard holomorphic (or anti-holomorphic) Whittaker function of weight $(2k_1,\cdots,2k_d)$ (or $(-2k_1,\cdots,-2k_d)$) defined in Section \ref{s3,s1}. The number $\widehat{\psi}(a_f)$  is called the $a_f$-th Fourier coefficient of $\psi$.

\begin{prop} \label{s3,s3,p1}
For any anti-holomorphic form $\psi\in\pi$ of level $t(\delta^{-1})U_0(Nc(\chi)^2)t(\delta)$
$$(\psi,\varphi^{\ast})=\widehat{\psi}(1)(\varphi^{\ast},\varphi^{\ast}).$$
\end{prop}
\begin{proof}
By newform theory such a $\psi$ has the form
$$\psi=c_0\varphi+\sum_{a|c(\chi)^2}c_a\varphi_a,$$
where $a$ is a nontrivial class modulo local units.  As $W_a(1)=0$ for nontrivial $a$, so
$\widehat{\psi}(1)=c_0\widehat{\varphi}(1)=c_0$, and 
$$(\psi,\varphi^{\ast})=(\widehat{\psi}(1)\varphi+\sum_{a|c(\chi)^2}c_a\varphi_a,\varphi^{\ast})=\widehat{\psi}(1)(\varphi^{\ast},\varphi^{\ast}).$$
 Here we have used the facts that    
$(\varphi^{\ast},\varphi_a)=0$
for non-unit $a|c(\pi_{\chi})$, and $(\varphi^{\ast},\varphi^{\ast})=(\varphi^{\ast},\varphi)$ (from the local definition).
\end{proof}

\begin{prop} \label{s3,s3,p2}
Let us retain the above notations, then 
\begin{align} \label{s3,s3,3}
L(1/2,\pi\times\pi_{\chi})=\frac{M\cdot L_{f}(1,\omega) 2^{|S|}}{\mu(ND)^+}\int_{Z(\BA)G(F)^+\backslash G(\BA)^+}\varphi^{\ast}(g)\varphi_{\chi}(g)E(1/2,g;\phi_2)dg,
\end{align}
where $L_{f}(1,\omega)$ denotes the finite part of the $L$-function, $D=c(\chi)^2c(\omega)$, $\mu(ND)^+=\mu(U_0(ND)^+)$, $S$ is the set of finite places dividing $c(\omega)$, and $M$ is given by
\begin{align} \label{M}
M=\prod_{v|\infty}2^{k_v+r_v+1}\prod_{v\in\Sigma_1}G_2(k_v-r_v)\prod_{v\in\Sigma_2}G_2(1+r_v-k_v).
\end{align}
\end{prop}
\begin{proof}
By Lemma \ref{s3,s2,l4} and newform theory we know that $\varphi^{\ast}=\varphi+\sum_{a|c(\chi)^2}c_a\varphi_a$, where the sum is over non-unit integral $a$ modulo the local units. Now (\ref{s3,s3,3}) follows from Proposition \ref{s3,s2,p4}, Proposition \ref{s3,s1,p1} and Lemma \ref{s3,s3,l1}. 
\end{proof}

\section{Theta correspondence and main formula} \label{s4}
In this section we first realize the terms $\varphi_{\chi}$ and $E(1/2,g;\phi_2)$ in  (\ref{s3,s3,p2}) as theta lifts. Then we apply a seesaw duality to express the the central value as a double torus integral. At last we use the Shimizu correspondence to determine the integrand explicitly. The main central value formula is therefore obtained.

\subsection{ Theta series and Eisenstein series} \label{s4,s1} 
Let $K$ as before be the quadratic CM extension of $F$. Let  $\phi_1=\prod_v\phi_{1,v}\in S(K_{\BA})$, where $\phi_{1,v}$ are  defined in Section \ref{s2,s2}. The theta lift  of $\chi$ from $GSO(K_{\BA})$ to $G(\BA)^+$   is defined by
\begin{align} \label{s4,s1,1}
\theta_{}(g,\chi;\phi_1)=\int_{SO(K)\backslash SO(\BA_K)}\theta_{}(t h,g;\phi_1)\chi(th)dt,
\end{align}
where $h$ is any element of $GSO(\BA_K)=\BA_K^{\times}$ such that $\N_K(h)=\det(g)$, and $\theta(h,g;\phi_1)$ is the theta kernel defined in  (\ref{s2,s1,2}) attached to $(K,\N_K)$. We
normalize the measure on $SO(\BA_K)=\BA_K^1$ such that $K_v^1\cap
U_{K_v}$  is of measure $1$ for any finite place $v$, and  $K_{v}^1=SO(2)$ is of measure
$1$ for any archimedean $v$. The total measure of $SO(K)\backslash SO(\BA_K)$ is $2L_f(1,\omega)$ (Waldspurger \cite{Waldspurger85} Section 1.5).

\begin{prop} \label{s4,s1,p1}
Let $\phi_2\in S(K_{v})$ be a Schwartz function in $S(K_v)$  defined by
\begin{equation*}
\phi_{1,v}(z)=\begin{cases}1_{\mathcal{O}_{K_v}}(z) &\text{ if $v$ is finite and $\chi_v$ is unramified,}\\
\chi_v(z)1_{\mathcal{O}_{K_v}^{\times}}(z) & \text{ if $\chi_v$ is ramified,}\\
2\bar{z}^{2r}e^{-2\pi (|z|^2)} &\text{ if $v$ is archimedean.}
\end{cases}
\end{equation*} 
(1) Let $ U_0(D)_v^{+}=U_0(D)\cap G(F_v)^+$ and let $\delta_v$  be the conductor of $\psi$ at a finite place $v$. For $k'\in t(\delta_v)^{-1} U_0(D)_v^+t(\delta_v)$ and  $g\in G(\BA)^+$ we have
$$\theta(\chi,gk';\phi_1)=\theta(\chi,g;\phi_1)$$
if $v$ is unramified in $K$, and
$$\theta(\chi,gk';\phi_1)=\omega(k')\theta(\chi,g;\phi_1)$$
if $v$ is ramified in $K$. Here in the last identity $\omega(k')=\omega(d)$ if $k'=t(\delta_v^{-1})\begin{pmatrix}a&b\\ c&d\end{pmatrix}t(\delta_v)$. Recall that $t(\delta)=\begin{pmatrix}\delta&0\\0&1\end{pmatrix}$.

At an archimedean $v$ the function $\theta(\chi,g;\phi_1)$ has weight $2r_v+1$. 

\noindent (2) For all  places $v$ not dividing $D$ and $g\in G(F_v)^+$, the local Whittaker function $W_{v}(\chi,g;\phi_1)$ of  $\theta_{}(\chi,g;\phi_1)$ equals $W_{\chi,v}(g)$ (defined in Section \ref{s3}). At $v| D$ we have $W_v(\chi,g;\phi_1)=W_{\chi,v}(g)$ for $g\in T_1(F_v)^+$.
\end{prop}

\begin{proof}
(1) The conclusion is already proved in \cite{Popa03} Proposition 4.1.2 for a finite $v$. The weight at an infinite $v$ is given in Proposition \ref{s2,s2,p1}.

(2)  By \cite{Popa03} Proposition 4.1.2 we  only need to  check $g=t(a)=\begin{pmatrix}a&0\\0&1\end{pmatrix}$
for $a>0$ at an archimedean place $v$. 

We quote a formula for the Whittaker function of $\theta(\chi,g;\phi)$ (see \cite{Shimizu72})
$$W_{v}(\chi,g;\phi_1)=\int_{\BC^1}L(h)r_1(g_1)\phi_1(\sigma^{-1})\chi(\sigma h)d\sigma,$$
where $h\in\BC$ is any number such that $\N(h)=\det(g)$. We will simply take $h=\sqrt{a}$.

As $\phi_1(z)=2(x-iy)^{2r}e^{-2\pi (x^2+y^2)}$, so we have
\begin{align*}
W_{v}\left(\chi,t(a);\phi_1\right)=&a^{1/2}\int_{\BC^1}\phi_1({a}t^{-1}h^{-1})\chi(th)dt\\
=&a^{1/2}\int_{\BR/2\pi\BZ}2(\overline{\sqrt{a}e^{-i\theta}})^{2r}e^{-2\pi a}e^{-2ri\theta}da\\
=&2a^{(2r+1)/2}e^{-2\pi a},
\end{align*}
which is exactly $W_{\chi,v}(t(a))$.
\end{proof}

Let $\theta(h,g;\phi_2)$ be the theta kernel attached to the quadratic space $(K,\Lambda\N_K)$ (Section \ref{s2,s1}). The theta lift of $1$ from $GSO(K_{\BA})$ to $G(\BA)^+$ is given by
$$I(g;\phi_2)=\theta(1,g;\phi_2)=\int_{SO(K)\backslash SO(K_{\BA})}\theta(th,g;\phi_2)dt,$$
where $h\in \BA_{K}^{\times}$ with $\N_K(h)=\det g$, and the measure on $SO(K)\backslash SO(K_{\BA})$ is the same as that used in defining the theta series $\theta(g,\chi;\phi_1)$. 

\begin{prop}[Siegel-Weil] \label{s4,s1,p3}
 For any $g\in G(\BA)^+$ one has
$$E(1/2,g;\phi_2)=L_f(1,\omega)^{-1} I(g;\phi_2).$$
\end{prop}  
\begin{proof}
For $g\in SL_2(\BA)$ and the Eisenstein series defined through the standard maximal compact subgroup,  the above formula is proved in Proposition 31 of \cite{Waldspurger80}.  At $s=1/2$ the standard Eisenstein series and our twisted one coincide as $|a(g)|^{2s-1}=|a_{\delta}(g)|^{2s-1}$ for $s=1/2$. 

The extension to $G(\BA)^+$ now follows from Theorem 4.2 of \cite{HK04}. 
\end{proof}

\begin{prop} \label{s4,s1,p4}
Let us retain the notations, then 
\begin{align} \label{s4,s1,2}
L(1/2,\pi\times\pi_{\chi})=\frac{2^{|S|}M} {\mu(ND)^+}\int_{Z(\BA)G(F)^+\backslash G(\BA)^+}\varphi^{\ast}(g)\theta(\chi,g;\phi_1)I(g;\phi_2)dg.
\end{align}
\end{prop}

\begin{proof}
The same argument as  \cite{Popa03} Proposition 4.3.1  shows the following identity for $\Re(s)\gg0$
\begin{align*}
\int_{Z(\BA)G(F)^+\backslash G(\BA)^+}\varphi'(g)\varphi_{\chi}(g)E_{}(s,g;\phi_2)dg=\int_{Z(\BA)G(F)^+\backslash G(\BA)^+}\varphi'(g)\theta(\chi,g;\phi_1)E_{}(s,g;\phi_2)dg,
\end{align*}
for any automorphic form $\varphi'$ in $\pi$. By  analytic continuation and Propositions \ref{s3,s3,p2}  
$$L(1/2,\pi\times\pi_{\chi})=\frac{2^{|S|}ML_f(1,\omega)} {\mu(ND)^+}\int_{Z(\BA)G(F)^+\backslash G(\BA)^+}\varphi^{\ast}(g)\theta(\chi,g;\phi_1)E_{}(1/2,g;\phi_2)dg.$$
Now Proposition \ref{s4,s1,p3} gives the desired formula. 
\end{proof}

\subsection{Seesaw identity} \label{s4,s2}
 We fix an embedding of $K$ into $B$ and $j\in B$ as in Section \ref{s2,s2}, such that $B=K+ Kj$ with $j^2=-\Lambda$. Therefore $(B,\N_B)=(K,\N_{K})\oplus (Kj,\N_K)=(K,\N)\oplus (K,\Lambda \N)$.

In Proposition \ref{s4,s1,p4} we find that the central value is an integral which involves two theta lifts $\theta(\chi,g;\phi_1)$ and $I(g;\phi_2)$.  To place these theta lifts together we use the following seesaw dual pair (\cite{Roberts98}):

\begin{equation} \label{s4,s2,g1}
\xymatrix{
R(G(\BA)^+\times G(\BA)^+)  \ar@{-} [dr]   & GSO(B_{\BA}) \ar@{-} [dl]   \\
G(\BA)^+ \ar [u]  & R(GSO(K_{\BA})\times GSO(K_{\BA}j)) \ar [u] \\
}
\end{equation}

Here the left vertical map is the diagonal embedding, and the right vertical one is the natural embedding given by regarding $(\mu,\nu)$ as the similitude
$$x+yj\mapsto \mu(x)+\nu(y)j\in B^{\times}_{\BA}.$$ 
Let  $r_1$, $r_{\Lambda}$ and  $r_B$  be the Weil representations defined on $R(GSO(K_{\BA})\times G(\BA)^+)$, $R(GSO(K_{\BA}j)\times G(\BA)^+)$, and on $R(GSO(B_{\BA})\times G(\BA)^+)$ respectively (Section \ref{s2,s1}).

We take  $\phi\in S(B_{\BA})=S(K_{\BA})\otimes S(K_{\BA}j)$ to be $\phi(x_1\oplus x_2j)=\phi_1(x_1)\phi_2(x_2)$, where $\phi_1$ and $\phi_2$ are given in Proposition \ref{s4,s1,p4}, so
$$r_B[(h_1,h_2),g]\phi(x_1\oplus x_2j)=r_1(h_1,g)\phi_1(x_1)r_{\Lambda}(h_2,g)\phi_2(x_2).$$
where $(h_1,h_2)\in R(GSO(K_{\BA})\times GSO(K_{\BA}))$ with similitude factors $\nu(h_1)=\nu(h_2)=\det(g)$. 
Consequently one   has a   decomposition for the theta kernels
\begin{equation} \label{s4,s2,1}
\theta_{B}((h_1,h_2),g;\phi)=\theta(h_1,g;\phi_1)\theta(h_2,g;\phi_2).
\end{equation}

Before going anywhere further we first prove a general seesaw identity for the above seesaw dual pair.   Let  $F_1$ and $F_2$ be two cuspidal forms on $Z(\BA)G(F)^+\backslash G(\BA)^+$ and $Z(\BA)H(F)\backslash H(\BA)$  respectively, here we write $H=R(GSO(K)\times GSO(Kj))$ for short. We define theta lifts
 $$\theta(F_1,h;\phi)=\int_{SL_2(F)\backslash SL_2(\BA)}\theta_B(h,g_1g;\phi)F_1(g_1g)dg_1,$$
$$\theta(F_2,g;\phi)=\int_{H_1(F)\backslash H_1(\BA)}\theta_B(h_1h,g;\phi)F_2(h_1h)dh_1,$$ where $\nu(h)=\det g$,

\begin{lem}[Seesaw identity]
  We have
\begin{equation} \label{s4,s2,2}
\int_{Z(\BA)G(F)^+\backslash G(\BA)^+}\theta(F_2,g;\phi)F_1(g)dg=\int_{Z(\BA)H(F)\backslash H(\BA)}\theta(F_1,h;\phi)F_2(h)dh.
\end{equation}
\end{lem}
\begin{proof}
Let $S$ be the compact group $F^+\BA_{\infty}^+\backslash\BA^+$, where $F^+=\N(K^{\times})$, $\BA^+=\N(\BA^{\times}_K)$ and $\BA_{\infty}^+=\N(K_{\infty}^{\times})$. Using the maps $\det$ of $G$ and $\nu$ of $H$ we have 
 $$1\to SL_2(F)\backslash SL_2(\BA)\to \BA_{\infty}^+G(F)^+\backslash G(\BA)^+\to S\to 1,$$and 
$$1\to H_1(F)\backslash H_1(\BA)\to \BA_{\infty}^+H(F)\backslash H(\BA)\to S\to1.$$ 
We also have
$$1\to \BA_{\infty}^+F^{\times}\backslash \BA^{\times}\to\BA_{\infty}^+G(F)^+\backslash G(\BA)^+\to Z(\BA)G(F)^+\backslash G(\BA)^+\to 1,$$
and
$$1\to \BA_{\infty}^+F^{\times}\backslash \BA^{\times}\to\BA_{\infty}^+H(F)\backslash H(\BA)\to Z(\BA)H(F)\backslash H(\BA)\to1.$$
Here the measure $ds$ on $S$ is the product of local measures on $K_v^+$ such that $\N_K(\mathcal{O}_{K,v}^{\times})$ has volume $1$. The measure on $\BA^+_{\infty}F^{\times}\backslash \BA^{\times}$ is the product of local measures on $F_v^{\times}$ such that $\mathcal{O}_v^{\times}$ has volume $1$. The measures on $SL_2(\BA)$ and $H_(F)\backslash H_1(\BA)$ are induced from the exact sequences. Since $H_1(F)\backslash H_1(\BA)=[SO(K)\backslash SO(K_{\BA})]^2$, the measure thus chosen on  $SO(K)\backslash SO(K_{\BA})$ is the same as the one in Section \ref{s4,s1}.  If we write $\mu=\mu(\BA^+_{\infty}F^{\times}\backslash \BA^{\times})$, then the left hand side of (\ref{s4,s2,2}) equals 
\begin{align*}
&\int_{Z(\BA)G(F)^+\backslash G(\BA)^+}\theta(F_2,g;\phi)F_1(g)dg\\
&=\frac1{\mu}\int_{\BA_{\infty}^+G(F)^+\backslash G(\BA)^+}\theta(F_2,g;\phi)F_1(g)dg\\
&=\frac1{\mu}\int_{S}\int_{SL_2(F)\backslash SL_2(\BA)}\theta(F_2,g_1s;\phi)F_1(g_1s)dg_1ds,
\end{align*}
and the right hand side of (\ref{s4,s2,2}) becomes
$$\frac1{\mu}\int_{S}\int_{H_1(F)\backslash H_1(\BA)}\theta(F_1,h_1s;\phi)F_2(h_1s)dh_1ds.$$
Now (\ref{s4,s2,2}) follows by Fubini theorem  and the definition of theta lifts.
\end{proof}
Now we take $F_1=\varphi^{\ast}$ and $F_2(h)=\chi(h)=\chi(t_1)$ for $h=(t_1,t_2)$, then by (\ref{s4,s2,1}) and (\ref{s4,s2,2}) the integral in (\ref{s4,s1,2}) becomes
\begin{align} \label{s4,s2,3}
\int_{Z(\BA)G(F)^+\backslash G(\BA)^+}\varphi^{\ast}(g)\theta(\chi,g;\phi_1)I(g;\phi_2)dg
=\int_{Z(\BA)H(F)\backslash H(\BA)}\theta(\varphi^{\ast},h;\phi)\chi(h)dh.
\end{align}

To obtain a better form of the integral we  use the following identification:
\begin{equation} \label{s4,s2,g2}
\xymatrix{
GSO(B_{\BA}) & B^{\times}_{\BA}\times B^{\times}_{\BA}/\BA^{\times} \ar [l] _{\sim} \\ H(\BA) \ar [u] & K^{\times}_{\BA}\times K^{\times}_{\BA}/\BA^{\times} \ar [l]  \ar [u]\\
}
\end{equation}
Here, the top isomorphism is given by $(g_1,g_2)b=g_1bg_2^{-1}$ for
$g_1,g_2\in B^{\times}$ and $b\in B$; the left vertical arrow is the one given in the above, and the right one is induced by the fixed embedding. So the bottom one sends $(t_1,t_2)\in K^{\times}_{\BA}\times K^{\times}_{\BA}$ to $(t_1t_2^{-1},t_1\bar{t}_2^{-1})=(h_1,h_2)\in R(GSO(_{\BA})\times GSO(K_{\BA}j))$. 

Under the above identification the right side of (\ref{s4,s2,3}) becomes

\begin{align} \label{s4,s2,4}
&\int_{Z(\BA)H(F)\backslash H(\BA)}\theta(\varphi^{\ast},h;\phi)\chi(h)dh
=\int_{(\BA^{\times}K^{\times}\backslash K_{\BA}^{\times})^2}\theta(\varphi^{\ast},(t_1,t_2);\phi)\chi(t_1t_2^{-1})dt_1dt_2.
\end{align} 

\begin{prop} \label{s4,s2,p1}
The central critical value of  $L(s,\pi\times\pi_{\chi})$ is given by
$$L(1/2,\pi_f\times\pi_{\chi})=\frac{2^{|S|}M} {\mu(ND)^+}\int_{(\BA^{\times}K^{\times}\backslash K^{\times}_{\BA})^2}\theta(\varphi^{\ast},(t_1,t_2);\phi)\chi(t_1t_2^{-1})dt_1dt_2,$$
where $M$ and $D$ are given in Proposition \ref{s3,s3,p1}.
\end{prop}

\subsection{Main formula} \label{s4,s3}
In this section we first show that the theta lift
$\theta(\varphi^{\ast},\sigma;\phi')$ on $\BA^{\times}\backslash B_{\BA}^{\times}\times B^{\times}_{\BA}$ decomposes as a product of two automorphic forms on $B_{\BA}^{\times}/B^{\times}\BA^{\times}$ (here and later on we always identify $GSO(B)$ and $B^{\times}\times B^{\times}/F^{\times}$). Then we use the decomposition to derive the main  central value formula.

\begin{lem} \label{s4,s3,l1}
The automorphic form $\theta(\varphi^{\ast},x,y;\phi)$, regarded as a form on $B^{\times}_{\BA}\times B^{\times}_{\BA}$, has the following level (or weight) structures:\\
1.   $\theta(\varphi^{\ast},xk_1, yk_2;\phi')=\chi(k_1^{-1}k_2)\theta(\varphi^{\ast},x,y;\phi')$ for $k_1,k_2\in\hat{R}^{\times}$,\\
2.  $\theta(\varphi^{\ast},xk_{\alpha},yk_{\beta};\phi)=e^{2ri(\alpha-\beta)}\theta(\varphi^{\ast},x,y;\phi)$ for $k_{\alpha}=e^{i\alpha}$. See Section \ref{s2,s2} for the definition of $\widehat{R}^{\times}$ and $\chi$ on it.
\end{lem}
\begin{proof}
The theta kernel is given by
$$\theta_B((x,y),g;\phi')=\sum_{b\in B(F)}r_B((x,y),g))\phi'(b)
=\sum_{b\in B(F)}|\N_B(x^{-1}y)|r_B(g_1)\phi'(x^{-1}by),$$
where $(x,y)\in {B^{\times}_{\BA}}^2$ such that $\nu(x,y)=\N_B(x y^{-1})=\det g$. Now the claims follow from Propositions \ref{s2,s3,p1} and \ref{s2,s2,p1}.
\end{proof}

\begin{prop} \label{s4,s3,p1}
We have 
$$\theta(\varphi^{\ast},x,y;\phi')=C\overline{{\varphi}^{B}(x)}\cdot{\varphi^{B}(y)},$$
where $C$ is  certain constant  to be determined later in Theorem \ref{thm}, $\varphi^B$ is an automorphic form in the automorphic representation $\pi^B$,  which is determined up-to a constant multiple  by the following level structures:\\
(1) $\varphi^B$ has weight $2r_v$ at an archimedean place $v$, \\
(2)   the action of  $k\in \widehat{R}^{\times}$ is given by
 $$\varphi^B(xk)=\chi(k)\varphi^B(x).$$
\end{prop}
\begin{proof}
In \cite{Shimizu72} it was shown that $\theta(\varphi,x,y;\phi)$ is in the product of ${\pi}^B\otimes\tilde{\pi}^B={\pi}^B\otimes{\pi}^B$, where $\tilde{\pi}^B$ is the contragredient of $\pi^B$. By Theorem 2.4.3 of \cite{Zhang01} level structures (1) and (2)  determine an  automorphic form $\varphi^B$ in $\pi^B$ uniquely (up-to a constant multiple). The statement is now clear by Lemma \ref{s4,s3,l1}.
\end{proof}

\begin{prop} \label{s4,s3,p3}
The central value is given by
\begin{align} \label{s4,s3,1}
L(1/2,\pi\times\pi_{\chi})=\frac{M}{\mu(N c(\chi)^2)^+}\left|\int_{\BA^{\times}K^{\times}\backslash \BA_K^{\times}}\varphi^B\chi^{-1}(t)dt\right|^2.
\end{align}

where  $M$ is given in (\ref{M}), and $C$, $\varphi^B$ are given in Proposition \ref{s4,s3,p1}.
\end{prop}
\begin{proof}
By Proposition \ref{s4,s3,p1} it suffices to show
\begin{align*}
&\int_{(\BA^{\times}K^{\times}\backslash K^{\times}_{\BA})^2}\theta(\varphi^{\ast},(t_1,t_2);\phi')\chi(t_1t_2^{-1})dt_1dt_2
\\ =& \frac{2^{|S|}}{\mu(c(\omega))^+}\int_{(\BA^{\times}K^{\times}\backslash K^{\times}_{\BA})^2}\theta(\varphi^{\ast},(t_1,t_2);\phi)\chi(t_1t_2^{-1})dt_1dt_2.
\end{align*}
But this is proved in Theorem 5.3.9 of \cite{Popa03}.  
\end{proof}

We now determine the constant $C$ in Proposition \ref{s4,s3,p3}.  The method used here  is inspired by \cite{Watson02}, but is simpler.

Let $\varphi'$ be the theta lift of ${\varphi^B (x)}\overline{\varphi^B(y)}$ from $GSO(B_{\BA})$ to $G(\BA)^+$ with respect to $\phi'\in S(B_{\BA})$ (here $G(\BA)^+$ denotes the matrices with determinants in $\N_B(B^{\times}_{\BA})$).
In other words (\cite{Shimizu72} or \cite{Watson02})
$$\varphi'(g)=\int_{\BA^{\times}B^{\times}\backslash
  B^{\times}_{\BA}}\int_{B^1\backslash
  B^1_{\BA}}\theta_B[(yx\sigma,y),g;\phi']{\varphi^B(yx\sigma)}\overline{\varphi^B(y)}dxdy,$$
with $\sigma\in B^{\times}_{\BA}$ such that $\N(\sigma)=\det(g)$.  
 The
measure on $GSO(B_{\BA})\cong PB^{\times}_{\BA}\times B^1_{\BA}$ is normalized such that the following adjoint identity holds
\begin{align} \label{s4,s5,-1}
({\varphi'},\overline{\varphi^{\ast}})=({\varphi^B}\overline{\varphi^B},\overline{\theta(\cdot,\varphi^{\ast},\phi')})_{GSO(B)}=C({\varphi^B}\overline{\varphi^B},{\varphi^B}\overline{\varphi^B}).
\end{align}
 
Shimizu \cite{Shimizu72}   showed that ${\varphi'}$ is  a cuspidal form in $\pi$ (after extending ${\varphi'}$ to $G(\BA)$ by left-invariance under $G(F)$).  By Proposition \ref{s2,s2,p1} and Proposition \ref{s2,s3,p1} the form $\varphi'$ has weight $(2k_1,\cdots,2k_d)$ and level $t(\delta)^{-1}U_0(N c(\chi)^2)t(\delta)$, so $\overline{\varphi}'$ is anti-holomorphic of weight $(-2k_1,\cdots,-2k_d)$. Proposition \ref{s3,s3,p1} implies
$$(\varphi',\overline{\varphi}^{\ast})=(\varphi^{\ast},\overline{\varphi}')=\overline{\widehat{\overline{\varphi}'}(1)}(\varphi^{\ast},\varphi^{\ast}).$$
The Whittaker function of $\overline{\varphi}'$ is related to that of $\varphi'$ by
$W_{\overline{\varphi}'}(g)=\overline{W_{\varphi'}(\epsilon g)}$,
where $\epsilon=t(-1)\in G(\BA)$, so $\widehat{\overline{\varphi}'}(1)=\overline{\widehat{\varphi}'(1)}$ and
\begin{align} \label{s4,s5,0}
(\varphi',\overline{\varphi}^{\ast})=\widehat{\varphi}'(1)(\varphi^{\ast},\varphi^{\ast}).
\end{align}

To compute the first Fourier coefficient of $\varphi'$ we note that the Whittaker function of ${\varphi'}$ is given by (\cite{Shimizu72}  \cite{Watson02})
\begin{align} \label{w}
{W_{\varphi'}}(g)=\int_{\BA^{\times}B^{\times}\backslash
  B^{\times}_{\BA}}\int_{B^1_{\BA}}L(x\sigma,1)r_B(g_1)\phi'(1){\varphi^B(yx\sigma)}dx\overline{\varphi^B(y)}dy.
\end{align}

\begin{prop} \label{s4,s5,p4}
Let
$$W(g,y)=\int_{B^1_{\BA}}L(x\sigma,1)r_B(g_1)\phi'(1){\varphi^B(yx\sigma)}dx,$$ then for $g=t(a_{\infty}1_f)$
 $$W(g,y)=\mu(\widehat{R}^1)M'W_{\infty}(t(a_{\infty}))\varphi^B(y),$$  
where $W_{\infty}$ is the standard Whittaker function of weight $(2k_1,\cdots,2k_d)$, $\mu(\widehat{R}^1)$ is the measure of $\widehat{R}^1\subset B^1_{\BA}$, and 
\begin{align} \label{M'}
M'=\prod_{v|\infty}2(4\pi)^{k_v-r_v-1}\frac{(k_v+r_v-1)!}{(2k_v-1)!}.
\end{align}
\end{prop}
\begin{proof}
At a finite place $v\not|c(\chi)$
\begin{align} \label{s4,s5,1}
&\int_{B^1_{F_v}}L(x,1)r_B(1)\phi'(1){\varphi^B(yx_v)}dx_v\\=&\int_{B^1_{F_v}}\phi'(x_v^{-1})\varphi^B(yx_v)dx_v
=\int_{R^1_v}\varphi^B(y x_v)dx_v
=\mu(R^1_v)\varphi^B(y)\nonumber
\end{align}
as $\varphi^B$ is invariant under $R^{\times}_v$ (here $R_v$ may be $\tilde{R}_v$). At a place $v|c(\chi)$
\begin{align} \label{s4,s5,2}
&\int_{B^1_{F_v}}L(x,1)r_B(1)\phi'(1){\varphi^B(yx_v)}dx_v\\=&\int_{B^1_{F_v}}\phi'(x_v^{-1})\varphi^B(yx_v)dx_v
=\int_{R^1_v}\chi(x_v)^{-1}\varphi^B(y x_v)dx_v
=\mu(R^1_v)\varphi^B(y) \nonumber
\end{align}
as $\varphi^B$ is $\chi$-isotypic under the action of $R^{\times}_v$.

We assume now $v$ is an archimedean place.  Because  ${\varphi'}$ has the lowest weight $2k_v$ at each archimedean place $v$   it must be holomorphic over $v$ and $W_{\varphi'}(t(a_v))=W_{\infty}(t(a_v))=0$ for $a_v<0$. Therefore we assume that $g=t(a_v)$ for $a_v>0$ from now on.

{\bf Definite case.}\,  Let $B_v=\BH$ be the Hamiltonian quaternion and take
$\sigma=a^{1/2}\in\BH$. The integral $W(g,y)$ is a Hecke operator on $\varphi^B$ at  place $v$. We will compute its eigenvalues  using the model of $\pi$  that is given by matrix coefficients $\{t^l_{rs}\}$ on $SU(2)$, where $l=k-1$ and $s$ varies (see Appendix \ref{s5,s2}). Because the weight of $\varphi^B$ at $v$ is $2r$ the vector $t^l_{rr}$ is  the one in this model that corresponds to $\varphi^B$.  
One has
\begin{align} \label{s4,s5,3}
\nonumber W(g,y)&=\int_{SU(2)}a\phi(\sigma^{\iota}x^{\iota})\varphi^B(yx\sigma)dx\\ \nonumber
&=t_{rr}^{k-1}(1)^{-1}\int_{SU(2)}a\cdot2(a^{1/2}{u})^{2r}p_{k-r-1}(4\pi a|v|^2)e^{-2\pi
  a(|u|^2+|v|^2)}t^{k-1}_{rr}(x)dx\cdot\varphi^B(y)\\ \nonumber
&=2a^{r+1}e^{-2\pi a}\int_{SU(2)}{u}^{2r}p_{k-r-1}(4\pi
a|v|^2)t^{k-1}_{rr}(x)dx\cdot\varphi^B(y)\\ 
&=2a^{r+1}e^{-2\pi
  a}\int_0^{\pi}\left(\frac{\cos\theta+1}2\right)^rp_{k-r-1}(2\pi
a(1-\cos\theta))P_{rr}^{k-1}(\cos\theta)\sin\theta d\theta\cdot \varphi^B(y)
\end{align}
Here in the above identity we have used ($l=1-k$)
$$P^{l}_{rr}(z)=\frac{(-1)^{l+r}}{2^l}\frac{1}{(l+r)!}(1+z)^r\frac{d^{l+r}}{dz^{l+r}}[(1+z)^{l-r}(1-z)^{l+r}],$$
and  the fact that
$\D t_{rr}^l(1)=P_{rr}^l(1)=1$.
The coefficient of the highest degree term  of $P^l_{rr}(z)$ is
$$\frac{(2l)!}{2^l(l+r)!(l-r)!}=\frac{(2k-2)!}{2^{k-1}(k+r-1)!(k-r-1)!},$$
while the  coefficient of the highest degree term (in $\cos\theta$) of $\left(\frac{\cos\theta+1}2\right)^rp_{k-r-1}(2\pi
a(1-\cos\theta))$ is given by
$$\frac{(2\pi a)^{k-r-1}}{2^{r}(k-r-1)!}.$$
Therefore, using the orthogonality of Legendre polynomials (\ref{orth2}) we can see that the integral in (\ref{s4,s5,3}) is  given by
$$2(4\pi)^{k-r-1}a^{k-r-1}\frac{(k+r-1)!}{(2k-1)!}.$$
Hence after putting $v$ back in we get
\begin{align} \label{s4,s5,4}
W_v(g,y)&=4\cdot4^{k_v-r_v-1}\pi^{k_v-r_v-1}a^{k_v}\frac{(k_v+r_v-1)!}{(2k_v-1)!}e^{-2\pi
  a}\varphi^B(y)\\ \nonumber
&=2(4\pi)^{k_v-r_v-1}\frac{(k_v+r_v-1)!}{(2k_v-1)!}W_{\infty}(t(a_v))\varphi^B(y).
\end{align}

{\bf Indefinite case.}\,
We  assume first that $ k_v>0$. We will occasionally drop the subscript $v$ in the following.

This time we  use the model generated by matrix elements $t_{rs}^{-k}$  on the group $SL_2(\BR)$. The vector in the model that corresponds to $\varphi^B$ is also $t^l_{rr}$ by weight consideration.  
Let $g=t(a)$ ($a>0$) and  $\sigma=\sqrt{a}$. The function $W(g,y)$ becomes:
\begin{align*}
W(g,y)=2a^{r+1}\int_{SL_2(\BR)}{u}^{2r}p_{r-k}(4\pi a|v|^2)e^{-2\pi a (|u|^2+|v|^2)}t^{-k}_{rr}(x) dx/t^{-k}_{rr}(1)\cdot\varphi^B(y).
\end{align*}
Using  Euler angles we obtain (Appendix \ref{s5,s2}):
$$W(g,y)=2a^{r+1}\int_{\BR^+}(\cosh\frac{t}2)^{2r}p_{r-k}(4\pi a\sinh^2(\frac{t}2))e^{-2a\pi\cosh t}{\frak P}^{-k}_{rr}(\cosh t)\sinh t dt,$$
and  we have
$${\frak P}^{-k}_{rr}(\cosh t)=(\cosh\frac{t}2)^{-2r}P_{r-k}^{(0,-2r)}(\cosh t).$$  

Note that $t^{-k}_{rr}(1)=1$, thus  the above integral becomes
\begin{align*}
W(g,y)=&2a^{r+1}\int_1^{\infty}p_{r-k}(2\pi a{(t-1)})e^{-2\pi a t}P_{r-k}^{(0,-2r)}dt\cdot \varphi^B(y)\\
=&2a^{r+1}\int_{0}^{\infty}p_{r-k}(2\pi a t)e^{-2\pi a}e^{-2\pi at}P_{r-k}^{(0,-2r)}(t+1)dt\cdot\varphi^B(y). 
\end{align*}
Recall that  the Jacobi polynomial $\D P^{(0,-2r)}_{r-k}(z)$ is given by:
$$\frac{\Gamma(r-k+1)}{(r-k)!}F(k-r,-k-r+1;1;\frac{1-z}2),$$
and $\D F(\alpha_1,\alpha_2;\beta;z)=\sum_{n=0}^{\infty}\frac{(\alpha_1)_n(\alpha_2)_n}{n!(\beta)_n}z^n.$ 
So we  get:
\begin{align*}W(g,h)&=2a^{r+1}e^{-2\pi a}\int_{\BR^+}p_{r-k}(2\pi at)\sum_{n=0}^{r-k}\frac{(k-r)_n(1-k-r)_n}{(n!)^2}\left(\frac{-t}2\right)^ne^{-2\pi a t}dt.\end{align*}
By the orthogonality relation of Laguerre polynomials  (\ref{orth1})
we can see  the integral equals:
\begin{align*}
&2a^{r+1}e^{-2\pi a}\frac{(k+r-1)\cdots (2k)}{(r-k)!}\int_{\BR_{+}}p_{r-k}(2\pi a t)\left(\frac{-t}2\right)^{r-k}e^{-2\pi a t}dt\\
=&2a^{r+1}e^{-2\pi a}\frac{(k+r-1)\cdots (2k)}{2^{r-k}(r-k)!}\cdot\frac{(r-k)!}{(2\pi a)^{r-k+1}}\\
=&2(4\pi)^{k-r-1}{(k+r-1)\cdots(2k)}2a^ke^{-2\pi a}.
\end{align*}
Putting   $v$ back and  we  obtain
\begin{align} \label{s4,s5,6}
W_v(g,y)=2(4\pi)^{k_v-r_v-1}\frac{(k_v+r_v-1)!}{(2k_v-1)!}W_{\infty}(t(a_v))\varphi^B(y).
\end{align}

 Proposition \ref{s4,s5,p4} now follows from (\ref{s4,s5,1}), (\ref{s4,s5,2}), (\ref{s4,s5,4}) and (\ref{s4,s5,6}).
\end{proof}

We now  come to the final central value formula. 
\begin{thm} \label{thm}
The central value of $L(s,\pi\times\pi_{\chi})$ is given by
$$L(1/2,\pi\times\pi_{\chi})=MM'\frac{\mu(\widehat{R}^1)}{\mu(Nc(\chi)^2)^+} \cdot \frac{(\varphi^{\ast},\varphi^{\ast})}{(\varphi^B,\varphi^B)}\left|\int_{\BA^{\times} K^{\times}\backslash \BA_{K}^{\times}}\varphi^B\chi^{-1}(t) dt\right|^2,$$
where $M$ is given by 
$$M=\prod_{v|\infty}2^{k_v+r_v+1}\prod_{v\in\Sigma_1}G_2(k_v-r_v)\prod_{v\in\Sigma_2}G_2(r_v-k_v+1),$$
$M'$ is given by
$$M'=\prod_{v|\infty}2(4\pi)^{k_v-r_v-1}\frac{(k_v+r_v-1)!}{(2k_v-1)!},$$
$\mu(\widehat{R}^1)$ is the measure of  $\widehat{R}^1=\widehat{R}^{\times}\cap B^1_{\BA_f}$ (Section \ref{s2,s2}) and $\mu(Nc(\chi)^2)^+$ is the measure of $U_0(Nc(\chi)^2)^+$. 
\end{thm}
\begin{proof}
From Proposition \ref{s4,s5,p4} 
$$\widehat{\varphi'}(1)=\mu(\widehat{R}^1)M'(\varphi^B,\varphi^B),$$
so by (\ref{s4,s5,-1}) and (\ref{s4,s5,0})
$$C=\widehat{\varphi}'(1)(\varphi^{\ast},\varphi^{\ast})/(\varphi^B,\varphi^B)^2=\mu(\widehat{R}^1)M'(\varphi^{\ast},\varphi^{\ast})/(\varphi^B,\varphi^B).$$
Now the formula is clear from Proposition \ref{s4,s3,p3}. 
\end{proof}

Note that the level of $\varphi^{\ast}$ is (the twist of) $N c(\chi)^2$, so is  the level of $\varphi^B$. In a subsequent paper we will lower the level to $N$, which is more suitable for applications.
\appendix

\section{Models of representations at infinity} \label{s5,s2}
In the next section we  need to compute  eigenvalues of certain archimedean Hecke operators which are defined as  integrals on $B_v^1$. Here $v$ is an archimedean place, and we will drop $v$ throughout this section. So $B^1$ is either $SU_2(\BC)$ or $SL_2(\BR)$ depending on whether $B$ is the quaternion algebra or the matrix algebra. To compute the eigenvalues we can  restrict the irreducible representation $\pi$ of $B^{\times}$ to $B^1$ (still denoted by $\pi$),  and calculate the eigenvalues on the restricted representation. This process will give the same eigenvalues. 

Now we describe a convenient model for the representation  $\pi$ of $B^1$. It realizes $\pi$  as a subspace of harmonic functions on the Lie group $B^1$. These functions are actually the matrix coefficients of $\pi$.  The facts  recorded here are well-known and  can be found, for instance, in Chapters 6 and 7 of \cite{Vilenkin91}.  We start with $B^1=SL_2(\BR).$

First, the group $SL_2(\BR)$ can be realized as a subgroup of $GL_2(\BC)$ and any  $g\in SL_2(\BR)$ is parametrized by three Euler angles
$$0\le\phi<2\pi, \,\, 0<t<\infty,\,\, -2\pi\le\psi<2\pi,$$
such that 
\begin{align} \label{a,1}
 g=\begin{pmatrix}e^{i\phi/2}&0\\0&e^{-i\phi/2}\end{pmatrix}
\begin{pmatrix}\cosh\frac{t}2&\sinh\frac{t}2\\\sinh\frac{t}2&\cosh\frac{t}2\end{pmatrix}\begin{pmatrix}e^{i\psi/2}&0\\0&e^{-i\psi/2}\end{pmatrix}.
\end{align} 
The measure of $SL_2(\BR)$ is given by 
$$\frac{1}{8\pi^2}\sinh t\,d\theta\, d\psi\, dt.$$
So the total measure of $SL_2(\BR)$ is one.
The torus $SO(2)$ is embedded in $SL_2(\BR)$ by mapping $e^{i\alpha}$ to $\begin{pmatrix}e^{i\alpha}&0\\0&e^{-i\alpha}\end{pmatrix}$.
 
Jacobi polynomials are defined by
$$P^{(\alpha,\beta)}_n(z)=\frac{\Gamma(n+\alpha+1)}{n!\Gamma(\alpha+1)}F(-n,n+\alpha+\beta+1;\alpha+1;\frac{1-z}2)$$
which is a polynomial of degree $n$ in $z$, here $\alpha\in\BZ_{\ge0}$, $n\in\BZ_{+}$ and $\D F(a,b;\alpha;z)=\sum_{n=0}^{\infty}\frac{(a)_n(b)_n}{n!(\alpha)_n}z^n$ is a hypergeometric series, where $(\alpha)_n=\alpha(\alpha+1)\cdots(\alpha+n-1)$.

Assume  $\pi$ is the discrete series of lowest weight $-2l$ with $l$ a negative integer.  Let  the  matrix elements of $\pi$ be  given by $t^l_{mn}(g)=(T^lv_m,v_n)$, where $(\, ,\,)$ is an invarinat Hermitian form and $v_m$ is a weight $m$ vector of norm 1.  So $t^l_{mn}(g)=e^{-i(m\phi+n\psi)}{\frak P}^l_{mn}(\cosh t)$ under the decomposition (\ref{a,1}), then
$${\frak P}^l_{mn}(\cosh t)=(\sinh \frac{t}2)^{m-n}(\cosh\frac{t}2)^{m+n}P^{(m-n,m+n)}_{l-m}(\cosh t),$$
where $n\le m\le l<0$, for other cases we use the symmetric relations. Another symmetric relation is
$${\frak P}^l_{mn}(\cosh t)={\frak P}^l_{-m,-n}(\cosh t).$$
If $\pi$ is the discrete series representation of $SL_2(\BR)$ of lowest weight $-2l$, then for every fixed $m\le l$   it has a model generated by the matrix coefficients $t^{l}_{mn}(g)$, where $n$ runs through all integer numbers less than or equal to $l$.

Similar facts hold  for $SU(2)$. The group $SU(2)$ can also be regarded as a subgroup of $GL_2(\BC)$ and every element $g\in SU(2)$  has an Euler angle parametrization
\begin{align}
g=\begin{pmatrix}e^{i\phi/2}&0\\0&e^{-i\phi/2}\end{pmatrix}\begin{pmatrix}\cos\frac{\theta}2& i\sin\frac{\theta}2\\ i\sin\frac{\theta}2& \cos\frac{\theta}2\end{pmatrix}\begin{pmatrix}e^{i\psi/2}&0\\ 0& e^{-i\psi/2}\end{pmatrix},
\end{align} 
where 
$$0\le\phi<2\pi,\,\, 0\le\theta\le\pi,\,\,-2\pi\le\psi<2\pi.$$
The measure is chosen such that the total measure of $SU(2)$ is $1$. We fix the same embedding of $SO(2)$ into $SU(2)$. The matrix coefficients of an irreducible representation $\pi$ of $SU(2)$ with dimension $2l+1$ have the form
$$t_{mn}^l(g)=i^{m-n}e^{-i(m\phi+n\psi)}P^l_{mn}(\cos\theta),$$
where $P^l_{mn}$ is a Legendre polynomial and $|m|, |n|\le l$. For a fixed $n$, the space generated by $t^l_{nm}$, where $m=-l,\cdots,l$, is a model of $\pi$ (under the right regular action).
We have the following Rodrigues formula for $P_{mn}^l(z)$:
\begin{align*}
P_{mn}^l(z)=&\frac{(-1)^{l-m}}{2^l}\left[\frac{(l+m)!}{(l-n)!(l+n)!(l-m)!}\right]^{1/2}\times\\
&\times (1+z)^{-(m+n)/2}(1-z)^{(n-m)/2}\frac{d^{l-m}}{dz^{l-m}}[(1-z)^{l-n}(1+z)^{l+n}].
\end{align*}
The Legendre polynomials have nice orthogonal properties:
\begin{align} \label{orth2}
\int_{-1}^{1}P^l_{mn}(z)P^r_{st}(z)dz=\frac{2}{2l+1}\delta_{lr}\delta_{ms}\delta_{nt},
\end{align} 
where $\delta_{ab}=1$ or $0$ depending on $a=b$ or not.

\bibliographystyle{amsplain}

\providecommand{\bysame}{\leavevmode\hbox
to3em{\hrulefill}\thinspace}

\bibliography{ref}

\end{document}